\documentclass{amsart}

\usepackage{amscd}
\usepackage{amssymb}
\usepackage{amsxtra}
\usepackage{longtable}
\usepackage{calc}
\usepackage{epsfig}
\usepackage{dbnsymb}

\allowdisplaybreaks[2]

\theoremstyle{plain}
\newtheorem{thm}{Theorem}
\newtheorem{cor}{Corollary}
\newtheorem{prop}{Proposition}
\newtheorem{lem}{Lemma}

\theoremstyle{definition}
\newtheorem{defn}{Definition}

\theoremstyle{remark}
\newtheorem{rmk}{Remark}
\newtheorem{rmks}{Remarks}
\newtheorem{exmp}{Example}
\newtheorem{exmps}{Examples}

\renewcommand{\epsilon}{\varepsilon}
\renewcommand{\kappa}{\varkappa}
\renewcommand{\theta}{\vartheta}

\newcommand{\IC}{\ensuremath{\mathbb C}}

\newcommand{\IN}{\ensuremath{\mathbb N}}
\newcommand{\IQ}{\ensuremath{\mathbb Q}}

\newcommand{\IZ}{\ensuremath{\mathbb Z}}

\newcommand{\ch}{\ensuremath{\mathrm{ch}}}

\newcommand{\End}{\ensuremath{\operatorname{End}}}
\newcommand{\Ext}{\ensuremath{\operatorname{Ext}}}

\newcommand{\Sym}{\ensuremath{\mathrm{S}}}

\newcommand{\Jet}{\ensuremath{\mathrm{J}}}

\newcommand{\varcard}[1]{\ensuremath{\operatorname{card}}}
\newcommand{\ord}[1]{\ensuremath{\operatorname{ord}}}
\newcommand{\set}[1]{\ensuremath{\left\{#1\right\}}}

\newcommand{\rk}{\ensuremath{\operatorname{rk}}}

\newcommand{\sltwo}{\ensuremath{\mathfrak{sl}_2}}

\newcommand{\Sp}{\ensuremath{\mathrm{Sp}}}

\newcommand{\orth}{\ensuremath{\bot}}

\newcommand{\scal}[1]{\ensuremath{\left<#1\right>}}

\newcommand{\HH}{\ensuremath{\mathrm H}}

\newcommand{\abs}[1]{\ensuremath{\left|#1\right|}}
\newcommand{\norm}[1]{\ensuremath{\left\Vert#1\right\Vert}}
\newcommand{\diff}{\ensuremath{\mathrm d}}

\newcommand{\eval}[1]{\ensuremath{\left.#1\right|}}

\DeclareMathOperator{\pr}{pr}
\DeclareMathOperator{\id}{id}

\newcommand{\OX}{\ensuremath{{\mathcal O_X}}}
\newcommand{\Tang}{\ensuremath{\mathcal T}}
\newcommand{\TangX}{\ensuremath{\Tang_X}}

\newcommand{\RW}{\ensuremath{\operatorname{RW}}}

\def\trivalent{\mathrm{t}}

\newcommand{\BB}{\ensuremath{\mathcal B}}
\newcommand{\hBB}{\ensuremath{\hat{\BB}}}
\newcommand{\BBo}{\ensuremath{\mbox{}^\trivalent\BB}}
\newcommand{\hBBo}{\ensuremath{\mbox{}^\trivalent\hBB}}
\newcommand{\hl}{\ensuremath{\hat\ell}}

\newcommand{\partitions}{\ensuremath{\mathrm{P}}}

\newcommand{\clos}[1]{\left<#1\right>}
\newcommand{\cclos}[1]{\left<\left<#1\right>\right>}
\newcommand{\nwheel}{\tilde w}

\newcommand{\Hilb}[2]{\ensuremath{{#2}^{\left[#1\right]}}}
\newcommand{\Kummer}[2]{\ensuremath{{#2}^{\left[\left[#1\right]\right]}}}

\def\N{\IN}
\def\Z{\IZ}
\def\Q{\IQ}

\def\C{\IC}

\title[Calculation of Rozansky-Witten invariants]
{Calculation of Rozansky-Witten invariants on the Hilbert schemes of
  points on a K3 surface and the generalised Kummer varieties}
\author[M.\ A.\ Nieper-Wi\ss kirchen]{Marc A.~Nieper-Wi\ss kirchen}
\thanks{The author is supported by the Deutsche Forschungsgemeinschaft.}
\thanks{This work has been written under hospitality of the Department of Pure
Mathematics and Mathematical Statistics of the University of Cambridge}
\address{University of Cambridge \\ Department of Pure Mathematics and
  Mathematical Statistics \\ Centre for Mathematical Sciences \\ Wilberforce
  Road \\ Cambridge CB3 0WB \\ United Kingdom}
\email{man35@dpmms.cam.ac.uk}
\date{\today}

\begin{document}

\begin{abstract}
  For any holomorphic symplectic manifold $(X, \sigma)$, a closed Jacobi
  diagram with $2k$ trivalent vertices gives rise to a Rozansky-Witten class
  $$\RW_{X, \sigma}(\Gamma) \in \HH^{2k}(X, \OX).$$ If $X$ is irreducible,
  this defines a number $\beta_\Gamma(X, \sigma)$ by $\RW_{X, \sigma}(\Gamma)
  = \beta_\Gamma(X, \sigma) [\bar\sigma]^k$.
  
  Let $(\Hilb n X, \Hilb n \sigma)$ be the Hilbert scheme of $n$ points on a
  K3 surface together with a symplectic form $\Hilb n \sigma$ such that
  $\int_{\Hilb n X} (\Hilb n \sigma \Hilb n{\bar\sigma})^n = n!$. Further, let
  $(\Kummer n A, \Kummer n \sigma)$ be the generalised Kummer variety of
  dimension $2n - 2$ together with a symplectic form $\Kummer n \sigma$ such
  that $\int_{\Kummer n A} (\Kummer n \sigma \Kummer n{\bar\sigma})^n = n!$.
  J.~Sawon conjectured in his doctoral thesis that for every connected Jacobi
  diagram, the functions $\beta_\Gamma(\Hilb n X, \Hilb n \sigma)$ and
  $\beta_\Gamma(\Kummer n A, \Kummer n \sigma)$ are linear in $n$.
  
  We prove that this conjecture is true for $\Gamma$ being a connected Jacobi
  diagram homologous to a polynomial of closed polywheels. We further show how
  this enables one to calculate all Rozansky-Witten invariants of $\Hilb n X$
  and $\Kummer n A$ for closed Jacobi diagrams that are homologous to a
  polynomial of closed polywheels.

  It seems to be unknown whether every Jacobi diagram is homologous to a
  polynomial of closed polywheels. If indeed the closed polywheels generate
  the whole graph homology space as an algebrea, our methods will thus enable
  us to compute \emph{all}
  Rozansky-Witten invariants for the Hilbert schemes and the generalised
  Kummer varieties using these methods.
  
  Also discussed in this article are the definitions of the various graph
  homology spaces, certain operators acting on these spaces and their
  relations, some general facts about holomorphic symplectic manifolds and
  facts about the special geometry of the Hilbert schemes of points on
  surfaces.
\end{abstract}

\maketitle

\tableofcontents

\section{Introduction}

A compact \emph{hyperk\"ahler manifold $(X, g)$} is a compact Riemannian
manifold whose holonomy is contained in $\Sp(n)$. An example of such a
manifold is the K3 surface together with a Ricci-flat K\"ahler metric (which
exists by S.~Yau's theorem~\cite{yau78}). In~\cite{rozansky97}, L.~Rozansky
and E.~Witten described how one can associate to every vertex-oriented
trivalent graph $\Gamma$ an invariant $b_\Gamma(X)$ to $X$, henceforth called
a \emph{Rozansky-Witten invariant of $X$ associated to $\Gamma$}. In fact,
this invariant only depends on the homology class of the graph, so the
invariants are already defined on the level of the graph homology space $\BB$
(see e.g.~\cite{barnatan95} and this paper for more information about graph
homology).

Every hyperk\"ahler manifold $(X, g)$ can be given the structure of a K\"ahler
manifold $X$ (which is, however, not uniquely defined) whose K\"ahler metric is
just given by $g$. $X$ happens to carry a holomorphic symplectic two-form
$\sigma \in \HH^0(X, \Omega_X^2)$, whereas we shall call $X$ a
\emph{holomorphic symplectic manifold}. Now M.~Kapranov showed
in~\cite{kapranov99} that one can in fact calculate $b_\Gamma(X)$ from $(X,
\sigma)$ by purely holomorphic methods. 

The basic idea is the following: We can identify the holomorphic tangent
bundle $\TangX$ of $X$ with its cotangent bundle $\Omega_X$ by means of
$\sigma$. Doing this, the Atiyah class $\alpha_X$ (see~\cite{kapranov99}) of
$X$ lies in $\HH^1(X, \Sym_3 \Tang_X)$. Now we place a copy of $\alpha_X$ at
each trivalent vertex of the graph, take the $\cup$-product of all these
copies (which gives us an element in $\HH^{2k}(X, (\Sym_3 \Tang_X)^{\otimes
  2k})$ if $2k$ is the number of trivalent vertices), and finally contract
$(\Sym_3 \Tang_X)^{\otimes 2k})$ along the edges of the graph by means of the
holomorphic symplectic form $\sigma$. Let us call the resulting element
$\RW_{X, \sigma}(\Gamma) \in \HH^{2k}(X, \OX)$.  In case $2k$ is the complex
dimension of $X$, we can integrate this element over $X$ after we have
multiplied it with $[\sigma]^{2k}$. This gives us more or less $b_\Gamma(X)$.
The orientation at the vertices of the graph is needed in the process to get a
number which is not only defined up to sign.

There are two main example series of holomorphic symplectic manifolds, the
Hilbert schemes $\Hilb n X$ of points on a K3 surface $X$ and the generalised
Kummer varieties $\Kummer n A$
(see~\cite{beauville83}). Besides two further manifolds constructed by
K.~O'Grady in~\cite{ogrady99} and~\cite{ogrady00}, these are the only known
examples of \emph{irreducible} holomorphic symplectic manifolds up to
deformation. 

Not much work was done on actual calculations of these invariants on the
example series. The first extensive calculations were carried out by J.~Sawon
in his doctoral thesis~\cite{sawon99}. All Chern numbers are in fact
Rozansky-Witten invariants associated to certain Jacobi diagrams, called
\emph{closed polywheels}. Let $\mathcal W$ be the subspace spanned by these
polywheels in $\BB$. All Rozansky-Witten invariants associated to graphs lying
in $\mathcal W$ can thus be calculated from the knowledge of the Chern numbers
(which are computable in the case of $\Hilb n X$ (\cite{lehn01}) or $\Kummer n
A$ (\cite{kummer}). However, from complex dimension four on, there are graph
homology classes that do not lie in $\mathcal W$.  J.~Sawon showed that for
some of these graphs the Rozansky-Witten invariants can still be calculated
from knowledge of the Chern numbers, which enables one to calculate all
Rozansky-Witten invariants up to dimension five. His calculations would work
for all irreducible holomorphic manifolds whose Chern numbers are known.

In this article, we will make use of the special geometry of $\Hilb n X$ and
$\Kummer n A$. Doing this, we are able to give a method with enables us to
calculate all Rozansky-Witten invariants for graphs homology classes that lie
in the \emph{algebra} $\mathcal C$ generated by \emph{closed polywheels} in
$\BB$. The closed polywheels form the subspace $\mathcal W$ of the algebra
$\BB$ of graph homology. This is really a proper subspace. However, $\mathcal
C$, the algebra generated by this subspace, is much larger, and, as far as the
author knows, it is unknown whether $\mathcal C = \BB$, i.e.~whether this work
enables us to calculate \emph{all} Rozanky-Witten invariants for the main
example series.

The idea to carry out this computations is the following: Let $(Y, \tau)$ be
any irreducible holomorphic symplectic manifold. Then $\HH^{2k}(Y, \mathcal
O_Y)$ is spanned by $[\bar\tau]^k$. Therefore, every graph $\Gamma$ with $2k$
trivalent vertices defines a number $\beta_\Gamma(Y, \tau)$ by $\RW_{Y,
  \tau}(\Gamma) = \beta_\Gamma(Y, \tau) [\bar\tau]^k$. J.~Sawon has already
discussed how knowledge of these numbers for connected graphs is enough to
deduce the values of all Rozansky-Witten invariants.

For the example series, let us fix holomorphic symplectic forms $\Hilb n
\sigma$, respective $\Kummer n \sigma$ with $\int_{\Hilb n X} (\Hilb n \sigma
\Hilb n {\bar\sigma})^n = n!$ respective $\int_{\Kummer n A} (\Kummer n \sigma
\Kummer n {\bar\sigma})^n = n!$.
J.~Sawon conjectured the following:
\begin{quote}
  The functions $\beta_{\Gamma}(\Hilb n X, \Hilb n \sigma)$ and
  $\beta_{\Gamma}(\Kummer n X, \Kummer n \sigma)$ are linear in $n$ for
  $\Gamma$ being a connected graph.
\end{quote}
The main result of this work is the proof of this conjecture for the class of
connected graphs lying in $\mathcal C$ (see Theorem~\ref{thm:betaconn}). We
further show how one can calculate these linear functions from the knowledge
of the Chern numbers and thus how to calculate all Rozansky-Witten invariants
for graphs in $\mathcal C$.

We should note that we don't make any use of the IHX relation in our
derivations, and so we could equally have worked on the level of Jacobi
diagrams.

Let us finally give a short description of each section. In
section~\ref{sec:prel} we collect some definitions and results which will be
used in later on. The next section is concerned with defining the algebra of
graph homology and certain operations on this space. We define~\emph{connected
  polywheels} and show how they are related with the usual closed polywheels
in graph homology. We further exhibit a natural $\mathfrak {sl}_2$-action on
an extended graph homology space. In section~\ref{sec:sympl}, we first look at
general holomorphic symplectic manifolds. Then we study the two example series
more deeply. Section~\ref{sec:rwinv} defines Rozansky-Witten invariants while
the last section is dedicated to the proof of our main theorem and
explicit calculations.

\section{Preliminaries}
\label{sec:prel}

\subsection{Some multilinear algebra}

Let $\mathcal T$ be a tensor category (commutative and with unit). For any
object $V$ in $\mathcal T$, we denote by $\Sym^k V$ the coinvariants of
$V^{\otimes k}$ with respect to the natural action of the symmetric group and
by $\Lambda^k V$ the coinvariants with respect to the alternating action.
Further, let us denote by $\Sym_k V$ and $\Lambda_k V$ the invariants of both
actions.

\begin{prop}
  \label{prop:lambda3}
  Let $I$ be a cyclicly ordered set of three elements. Let $V$ be an object in
  $\mathcal T$. Then there exists a unique map $\Lambda_3 V \to V^{\otimes I}$
  such that for every bijection $\phi: \set{1, 2, 3} \to I$ respecting the
  canonical cyclic ordering of $\set{1, 2, 3}$ and the given cyclic ordering
  of $I$ the following diagram 
  \begin{gather}
    \begin{CD}
      \Lambda_3 V @= \Lambda_3 V \\
      @VVV @VVV \\
      V^{\otimes 3} @>>{\phi_*}> V^{\otimes I}
    \end{CD}
  \end{gather}
  commutes, where the map $\phi_*$ is the canonical one induced by $\phi$.
\end{prop}

\begin{proof}
  Let $\phi, \phi': \set{1, 2, 3} \to I$ be two bijections respecting the
  cyclic ordering. Then there exist an even permutation $\alpha \in \mathfrak
  A_3$ such that the lower square of the following diagram commutes:
  \begin{gather*}
    \begin{CD}
      \Lambda_3 V @= \Lambda_3 V \\
      @VVV @VVV \\
      V^{\otimes 3} @>>{\alpha_*}> V^{\otimes 3} \\
      @V{\phi}VV @VV{\phi'}V \\
      V^{\otimes I} @= V^{\otimes I}.
    \end{CD}
  \end{gather*}
  We have to show that the outer rectangle commutes. For this it suffices to
  show that the upper square commutes. In fact, since $\alpha$ is an even
  permutation, every element of $\Lambda_3 V$ is by definition invariant under
  $\alpha_*$.
\end{proof}

\subsection{Partitions}

A partition $\lambda$ of a non-negative integer $n \in \N_0$ is a sequence
$\lambda_1, \lambda_2, \dots$ of non-negative integers such that
\begin{gather}
  \norm \lambda := \sum_{i = 1}^\infty i \lambda_i = n.
\end{gather}
Therefore almost all $\lambda_i$ have to vanish.
In the literature, $\lambda$ is often notated by $1^{\lambda_1} 2^{\lambda_2}
\dots$. The set of all partitions of $n$ is denoted by $\partitions(n)$. The
union of all $\partitions(n)$ is denoted by $\partitions := \bigcup_{n =
  0}^\infty \partitions(n)$.
For every partition $\lambda \in \partitions$, we set
\begin{gather}
  \abs\lambda := \sum_{i = 1}^\infty \lambda_i
\end{gather}
and
\begin{gather}
  \lambda! := \prod_{i = 1}^\infty \lambda_i!.
\end{gather}

Let $a_1, a_2, \dots$ be any sequence of elements of a commutative unitary
ring. We set
\begin{gather}
  a_\lambda := \prod_{i = 1}^\infty a_i^{\lambda_i}
\end{gather}
for any partition $\lambda \in \partitions$.

With these definitions, we can formulate the following proposition in a nice
way:
\begin{prop}
  \label{prop:exppart}
  In $\Q[[a_1, a_2, \dots]]$ we have
  \begin{gather}
    \exp\left(\sum_{i = 1}^\infty a_i\right)
    = \sum_{\lambda \in \partitions} \frac{a_\lambda}{\lambda!}.
  \end{gather}
\end{prop} 

\begin{proof}
  We calculate
  \begin{gather}
    \exp\left(\sum_{i = 1}^\infty a_i\right)
    = \sum_{n = 0}^\infty \frac 1 {n!} \left(\sum_{i = 1}^\infty a_i\right)^n
    = \sum_{n = 0}^\infty \frac 1 {n!} \sum_{\lambda \in \partitions, \abs
      \lambda = n} n! \prod_{i = 1}^\infty \frac{a_i^{\lambda_i}}{\lambda_i!}
    = \sum_{\lambda \in \partitions} \frac{a_\lambda}{\lambda!}.
  \end{gather}
\end{proof}

If we set
\begin{gather}
  \frac{\partial}{\partial a_\lambda} := \prod_{i = 1}^\infty
  \eval{\frac{\partial^{\lambda_i}} {\partial a_i^{\lambda_i}}}_{a_i = 0},
\end{gather}
we have due to Proposition~\ref{prop:exppart}:
\begin{prop}
  In $\Q[[s_1, s_2, \dots]][a_1, a_2, \dots]$ we have
  \begin{gather}
    \frac{\partial}{\partial a_\lambda} \exp\left(\sum_{i = 1}^\infty a_i
      s_i\right) = s_\lambda.
  \end{gather}
\end{prop}

\subsection{A lemma from umbral calculus}

\begin{lem}
  \label{lem:sheffer}
  Let $R$ be any $\Q$-algebra (commutative and with unit) and $A(t) \in R[[t]]$
  and $B(t) \in tR[[t]]$ be two power series. Let the polynomial sequences
  $(p_n(x))$ and $(s_n(x))$ be defined by
  \begin{align}
    \label{eq:pdef}
    \sum_{k = 0}^\infty p_k(x) \frac{t^k}{k!} & = \exp(x B(t))
    \\
    \intertext{and}
    \label{eq:sdef}
    \sum_{k = 0}^\infty s_k(x) \frac{t^k}{k!} & = A(t) \exp(x B(t)).
  \end{align}
  Let $W_B(t) \in tR[[t]]$ be defined by $W_B(t \exp(B(t))) = t$.
  Then we have
  \begin{align}
    \label{eq:pstat}
    \sum_{k = 0}^\infty \frac{x p_k(x - k)}{(x - k)} \frac{t^k}{k!}
    & = \exp(x
    B(W_B(t))) \\
    \intertext{and}
    \label{eq:sstat}
    \sum_{k = 0}^\infty \frac{s_k(x - k)}{k!}t^k &
    = \frac{A(W_B(t))}{1 + W_B(t) B'(W_B(t))} \exp(x B(W_B(t))).
  \end{align}
\end{lem}

\begin{proof}
  It suffices to prove the result for the field $R = \Q(a_0, a_1, \dots, b_1,
  b_2, \dots)$ and $A(t) = \sum_{k = 0}^\infty a_k t^k$ and $B(t) = \sum_{k =
    1}^\infty b_k t^k$. 
  
  So let us assume this special case for the rest of the proof. Let us denote
  by $f(t)$ the compositional inverse of $B(t)$, i.e.~$f(B(t)) = t$. We set
  $g(t) := A^{-1}(f(t))$. For the following we will make use of the
  terminology and the statements in~\cite{roman84}. Using this
  terminology,~\eqref{eq:pdef} states that $(p_n(x))$ is the associated
  sequence to $f(t)$ and~\eqref{eq:sdef} states that $(s_n(x))$ is the Sheffer
  sequence to the pair $(g(t), f(t))$ (see Theorem~2.3.4 in~\cite{roman84}).

  Theorem~3.8.3 in~\cite{roman84} tells us that $(s_n(x-n))$ is the Sheffer
  sequence to the pair $(\tilde g(t), \tilde f(t))$ with
  \begin{align*}
    \tilde g(t) = g(t) (1 + f(t)/f'(t)) \\
    \intertext{and}
    \tilde f(t) = f(t) \exp(t).
  \end{align*}
  The compositional inverse of $\tilde f(t)$ is given by $\tilde B(t) :=
  B(W_B(t))$:
  \begin{gather*}
    B(W_B(\tilde f(t))) = B(W_B(f(t) \exp(t))) = B(W_B(f(t) \exp(B(f(t)))))
    = B(f(t)) = t.
  \end{gather*}
  Further, we have
  \begin{multline*}
    \tilde A(t) := \tilde g^{-1}(\tilde B(t))
    \\
    = (g(B(t)) (1 + f(B(t))/f'(B(t))))^{-1} \circ W_B(t)
    = \frac{A(t)}{1 + t B'(t)} \circ W_B(t),
  \end{multline*}
  which proves~\eqref{eq:sstat} again due to Theorem~2.3.4 in~\cite{roman84}.
  
  It remains to prove~\eqref{eq:pstat}, i.e.~that $(\frac {x p_n(x)}{x - n})$
  is the associated sequence to $\tilde f(t)$. We already know that $(p_n(x -
  n))$ is the Sheffer sequence to the pair $(1 + f(t)/f'(t), \tilde f(t))$. By
  Theorem~2.3.6 of~\cite{roman84} it follows that the associated sequence to
  $\tilde f(t)$ is given by $(1 + f(\diff/\diff x)/f'(\diff/\diff x)) p_n(x -
  n)$. By Theorem~2.3.7 and Corollary~3.6.6 in~\cite{roman84}, we have
  \begin{multline*}
    \left(1 + \frac{f(\diff/\diff x)}{f'(\diff/\diff x)}\right) p_n(x - n)
    = p_n(x - n) + \frac 1 {f'(\diff/\diff x)} n p_{n - 1}(x - n)
    \\
    = p_n(x - n) + \frac {n p_n(x - n)} {x - n}
    = \frac {x p_n(x - n)}{x - n},
  \end{multline*}
  which proves the rest of the lemma.
\end{proof}

\section{Graph homology}

This section is concerned with the space of graph homology classes of
unitrivalent graphs. A very detailed discussion of this space and other graph
homology spaces can be found in~\cite{barnatan95}. Further aspects of graph
homology can be found in~\cite{thurston00}, and, with respect to
Rozansky-Witten invariant, in~\cite{hitchin01}.

\subsection{The graph homology space}

In this article, \emph{graph} means a collection of vertices connected by
edges, i.e.\ every edge connects two vertices.  We want to call a half-edge
(i.e.\ an edge together with an adjacent vertex) of a graph a \emph{flag}. So,
every edge consists of exactly two flags. Every flag belongs to exactly one
vertex of the graph. On the other hand, a vertex is given by the set of its
flags. It is called \emph{univalent} if there is only one flag belonging to
it, and it is called \emph{trivalent} if there are exactly three flags
belonging to it. We shall identify edges and vertices with the set of their
flags. We shall also call univalent vertices \emph{legs}.
\begin{figure}[h]
  \begin{center}
    \epsfig{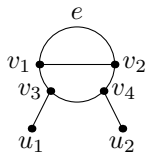}
  \end{center}
  \caption{This Jacobi diagram has four trivalent
    vertices $v_1, \dots, v_4$, and two univalent vertices $u_1$ and $u_2$,
    and $e$ is one of its $7$ edges.}
\end{figure}
A graph is called \emph{vertex-oriented} if, for every vertex, a
cyclic ordering of its flags is fixed.

\begin{defn}
  A \emph{Jacobi diagram} is a vertex-oriented graph with only uni-
  and trivalent vertices. A \emph{connected Jacobi diagram} is a
  Jacobi diagram which is connected as a graph. A \emph{trivalent
    Jacobi diagram} is a Jacobi diagram with no univalent vertices.
  
  We define the \emph{degree of a Jacobi diagram} to be the number of
  its vertices. It is always an even number.

  We identify two graphs if they are isomorphic as vertex-oriented graphs in
  the obvious sense.
\end{defn}
 
\begin{exmp}
  The empty graph is a Jacobi diagram, denoted by $1$.
  The unique Jacobi diagram consisting of two univalent vertices (which
  are connected by an edge) is denoted by $\ell$.
  \begin{figure}[h]
    \begin{center}
      \epsfig{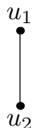}
    \end{center}
    \caption{The Jacobi diagram $\ell$ with its two univalent
      vertices $u_1$ and $u_2$.}
  \end{figure}
\end{exmp}

\begin{rmk}
  There are different names in the literature for what we call a ``Jacobi
  diagram'', e.g.\ unitrivalent graphs, chord diagrams, Chinese characters,
  Feynman diagrams. The name chosen here is also used by D.~Thurston
  in~\cite{thurston00}. The name comes from the fact that the IHX relation in
  graph homology defined later is essentially the well-known Jacobi identity
  for Lie algebras.
  
  With our definition of the degree of a Jacobi diagram, the algebra of graph
  homology defined later will be commutative in the graded sense. Further, the
  map $\RW$ that will associate to each Jacobi diagram a Rozansky-Witten class
  will respect this grading. But note that often the degree is defined to be
  \emph{half} of the number of vertices, which still is an integer.
\end{rmk}

We can always draw a Jacobi diagram in a planar drawing so that it
looks like a planar graph with vertices of valence $1$, $3$ or
$4$. Each $4$-valent vertex has to be interpreted as a crossing of two
non-connected edges of the drawn graph and not as one of its vertices.
Further, we want the counter-clockwise ordering of the flags at each
trivalent vertex in the drawing to be the same as the given cyclic
ordering.
\begin{figure}[h]
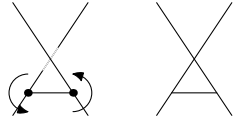

  \begin{center}
    \epsfig{file=graphs.3} \qquad \epsfig{file=graphs.4}
  \end{center}
  \caption{These two graphs depict the same one.}
\end{figure}
 
In drawn Jacobi diagrams, we also use a notation like $\cdots\overset n
-\cdots$ for a part of a graph which looks like a long line with $n$
univalent vertices (``legs'') attached to it, for example
$\ldots\bot\!\bot\!\bot\ldots$ for $n = 3$. The position of
$n$ indicates the placement of the legs relative to the ``long line''.

\begin{defn}
  Let $\mathcal T$ be any tensor category (commutative and with unit). Every
  Jacobi diagram $\Gamma$ with $k$ trivalent and $l$ univalent vertices
  induces a natural transformation $\Psi^\Gamma$ between the functors
  \begin{align}
    \mathcal T \to \mathcal T, & V \mapsto \Sym_k \Lambda_3 V \otimes
    \Sym_l V \\
    \intertext{and}
    \mathcal T \to \mathcal T, & V \mapsto \Sym^e \Sym^2 V,
  \end{align}
  where $e := \frac{3 k + l}{2}$ which is given by
  \begin{multline}
    \Psi^\Gamma: \Sym_k \Lambda_3 V \otimes \Sym_l V
    \overset{(1)} \to
    \bigotimes_{t \in T} \Lambda_3 V \otimes \bigotimes_{f \in U} V
    \overset{(2)} \to
    \bigotimes_{t \in T} \bigotimes_{f \in t} V
    \otimes \bigotimes_{f \in U} V
    \\
    \overset{(3)} \to
    \bigotimes_{f \in F} V
    \overset{(4)} \to
    \bigotimes_{e \in E}
    \bigotimes_{f \in e} V
    \overset{(5)} \to
    \Sym^e \Lambda^2 V,
  \end{multline}
  where $T$ is the set of the trivalent vertices, $U$ the set of the univalent
  vertices, $F$ the set of flags, and $E$ the set of edges of
  $\Gamma$. Further,
  \begin{enumerate}
  \item[(1)] is given by the natural inclusions of the invariants in the
    tensor products,
  \item[(2)] is given by the canonical maps (see
  Proposition~\ref{prop:lambda3} and recall that the sets $t$ are cyclicly
  ordered),
  \item[(3)] is given by the associativity of the tensor product,
  \item[(4)] is given again by the associativity of the tensor product, and
  finally
  \item[(5)] is given by the canonical projections onto the coinvariants.
  \end{enumerate}
\end{defn}

\begin{defn}
  We define $\BB$ to be the $\Q$-vector space spanned by all Jacobi
  diagrams modulo the IHX relation
  \begin{gather}
    \IGraph = \HGraph - \XGraph
  \end{gather}
  and the anti-symmetry (AS) relation
  \begin{gather}
    \YGraph + \TwistedY = 0,
  \end{gather}
  which can be applied anywhere within a diagram. (For this definition see
  also~\cite{barnatan95} and~\cite{thurston00}.) Two Jacobi diagrams are said
  to be \emph{homologous} if they are in the same class modulo the IHX and AS
  relation.

  Furthermore, let $\BB'$ be the subspace of $\BB$ spanned by all Jacobi
  diagrams not containing $\ell$ as a component, and let $\BBo$ be the
  subspace of $\BB'$ spanned by all trivalent Jacobi diagrams. All
  these are graded and double-graded. The grading is induced by the
  degree of Jacobi diagrams, the double-grading by the number of
  univalent and trivalent vertices.

  The completion of $\BB$ (resp.\ $\BB'$, resp.\ $\BBo$) with respect
  to the grading will be denoted by $\hBB$ (resp.\ $\hBB'$, resp.\ $\hBBo$).
  
  We define $\BB_{k, l}$ to be the subspace of $\hBB$ generated by graphs with
  $k$ trivalent and $l$ univalent vertices. $\BB'_{k, l}$ and $\BBo_k :=
  \BBo_{k, 0}$ are defined similarly.
\end{defn}

All these spaces are called \emph{graph homology spaces} and their
elements are called \emph{graph homology classes} or \emph{graphs} for
short. 

\begin{rmk}
  The subspaces $\BB_k$ of $\hBB$ spanned by the Jacobi diagrams of degree $k$
  are always of finite dimension. The subspace $\BB_0$ is one-dimensional and
  spanned by the graph homology class $1$ of the empty diagram $1$.
\end{rmk}

\begin{rmk}
  We have $\hBB = \prod_{k, l \geq 0} B_{k, l}$. Further, $\hBB =
  \hBB'[[\ell]]$.
  Due to the AS relation, the spaces $\BB'_{k, l}$ are zero for $l >
  k$. Therefore, $\hBB' = \prod_{k = 0}^\infty \bigoplus_{l = 0}^k \BB'_{k,
  l}$.
\end{rmk}

\begin{exmp}
  If $\gamma$ is a graph which has a part looking like
  $\cdots\overset n -\cdots$, it will become $(-1)^n \gamma$ if we
  substitute the part $\cdots\overset n -\cdots$ by $\cdots\underset n
  -\cdots$ due to the anti-symmetry relation.
\end{exmp}

\subsection{Operations with graphs}

\begin{defn}
  Disjoint union of Jacobi diagrams induces a bilinear map
  \begin{gather}
    \hBB \times \hBB \to \hBB, (\gamma, \gamma') \mapsto \gamma \cup
    \gamma'.
  \end{gather}
\end{defn}

By mapping $1 \in \Q$ to $1 \in \hBB$, the space $\hBB$ becomes a
graded $\Q$-algebra, which has no components in odd degrees. Often, we
omit the product sign ``$\cup$''. $\BB$, $\BB'$, $\BBo$, and so on
are subalgebras.

\begin{defn}
  Let $k \in \N$. We call the graph homology class of the Jacobi
  diagram $\overset{2k}{\bigcirc}$ the \emph{$2k$-wheel $w_{2k}$},
  i.e.\ $w_2 = \twowheel$, $w_4 = \fourwheel$, and so on. It has $2k$
  univalent and $2k$ trivalent vertices. The expression $w_0$ will be given a
  meaning later, see section~\ref{sec:sltwo}.
\end{defn}

\begin{rmk}
  The wheels $w_k$ with $k$ odd vanish in $\hBB$ due to the AS relation.
\end{rmk}

Let $\Gamma$ be a Jacobi diagram and $u,
u'$ be two different univalent vertices of $\Gamma$. These two should not
be the two vertices of a component $\ell$ of $\Gamma$. Let $v$
(resp.\ $v'$) be the vertex $u$ (resp.\ $u'$) is attached to. The process of
\emph{gluing the vertices $u$ and $u'$} means to remove $u$ and $u'$
together with the edges connecting them to $v$ resp.\ $v'$ and to add
a new edge between $v$ and $v'$. Thus, we arrive at a new graph
$\Gamma/(u, u')$, whose number of trivalent vertices is the number of
trivalent vertices of $\Gamma$ and whose number of univalent vertices
is the number of univalent vertices of $\Gamma$ minus two. To make it
a Jacobi diagram we define the cyclic orientation of the flags at $v$
(resp.\ $v'$) to be the cyclic orientation of the flags at $v$ (resp.\ 
$v'$) in $\Gamma$ with the flag belonging to the edge connecting $v$
(resp.\ $v'$) with $u$ (resp.\ $u'$) replaced by the flag belonging to
the added edge.
\begin{figure}[h]
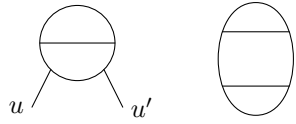

  \begin{center}
    \epsfig{file=graphs.6} \qquad \epsfig{file=graphs.7}
  \end{center}
  \caption{Gluing the two
    univalent vertices $u$ and $u'$ of the left graph produces the right one,
    denoted by $\ThetaGraph_2$.}
\end{figure}
For example, gluing the two univalent vertices of $w_2$ leads to the
graph $\ThetaGraph$.

If $\pi = \set{\set{u_1, u'_1}, \dots, \set{u_k, u'_k}}$ is a set of
two-element sets of legs that are pairwise disjoint and such that each pair
$u_k, u'_k$ fulfills the assumptions of the previous construction, we set
\begin{gather}
  \Gamma/\pi := \Gamma/(u_1, u'_1)/\dots/(u_k, u'_k).
\end{gather}

Of course, the process of gluing two univalent vertices given above
does not work if $u$ and $u'$ are the two univalent vertices of
$\ell$, thus our assumption on $\Gamma$.

\begin{defn}
  Let $\Gamma, \Gamma'$ be two Jacobi diagrams, at least one of them
  without $\ell$ as a component and $U = \set{u_1, \dots, u_n}$
  resp.\ $U'$ the sets of their univalent vertices.  We define
  \begin{gather}
    \hat\Gamma(\Gamma') := \sum_{\substack{f: U \hookrightarrow U' \\
    \text{injective}}} (\Gamma \cup \Gamma')/(u_1, f(u_1))/\dots/(u_n,
    f(u_n)),
  \end{gather}
  viewed as an element in $\hBB$.

  This induces for every $\gamma \in \hBB$ a $\hBBo$-linear map
  \begin{gather}
    \hat\gamma: \hBB' \to \hBB', \gamma' \mapsto \hat\gamma(\gamma').
  \end{gather}
\end{defn}

\begin{exmp}
  Set $\partial := \frac 1 2 \hl$.
  It is is an endomorphism of $\hBB'$ of degree $-2$. For example,
  $\partial \twowheel = \ThetaGraph$.
  By setting
  \begin{gather}
    \partial(\gamma, \gamma') := \partial(\gamma\cup\gamma') -
    \partial(\gamma)\cup\gamma' - \gamma\cup\partial(\gamma')
  \end{gather}
  for $\gamma, \gamma' \in \hBB'$,
  we have the following formula for all $\gamma \in
  \hBB'$:
  \begin{gather}
    \partial(\gamma^n) = \binom n 1 \partial(\gamma)\gamma^{n - 1} +
    \binom n 2 \partial(\gamma, \gamma) \gamma^{n - 2}.
  \end{gather}
  This shows that $\partial$ is a differential operator of order two acting on
  $\hBB'$.

  Acting by $\partial$ on a Jacobi diagram means to glue two of its
  univalent vertices in all possible ways, acting by $\partial(\cdot,
  \cdot)$ on two Jacobi diagrams means to connect them by gluing a
  univalent vertex of the first with a univalent vertex of the second
  in all possible ways.
\end{exmp}

\begin{defn}
  Let $\Gamma, \Gamma'$ be two Jacobi diagrams, at least one of them
  without $\ell$ as a component, and $U = \set{u_1, \dots, u_n}$
  resp.\ $U'$ the sets of their univalent vertices.  We define
  \begin{gather}
    \left<\Gamma, \Gamma'\right>
    := \sum_{\substack{f: U \to U' \\
        \text{bijective}}} (\Gamma \cup \Gamma')/(u_1, f(u_1))/\dots/(u_n,
    f(u_n)),
  \end{gather}
  viewed as an element in $\hBBo$.

  This induces a $\hBBo$-bilinear map
  \begin{gather}
    \left<\cdot, \cdot\right>: \hBB' \times \hBB \to \hBBo,
  \end{gather}
  which is symmetric on $\hBB' \times \hBB'$.
\end{defn}

Note that $\left<\Gamma, \Gamma'\right>$ is zero unless $\Gamma$ and
$\Gamma'$ have equal numbers of univalent vertices. In this case, the
expression is the sum over all possibilities to glue the univalent
vertices of $\Gamma$ with univalent vertices of $\Gamma'$.

Note that $\left<\Gamma, \Gamma'\right>$ is zero unless $\Gamma$ and
$\Gamma'$ have equal numbers of univalent vertices. In this case, the
expression is the sum over all possibilities to glue the univalent
vertices of $\Gamma$ with univalent vertices of $\Gamma'$.

\begin{prop}
  \label{prop:ellandpartial}
  The map $\left<1, \cdot\right>: \hBB \to \hBBo$ is the canonical
  projection map, i.e.\ it removes all non-trivalent components from a
  graph. Furthermore, for $\gamma \in \hBB'$ and $\gamma' \in \hBB$, we
  have
  \begin{gather}
    \label{equ:ellandpartial}
    \left<\gamma, \frac\ell 2 \gamma'\right> = \left<\partial \gamma,
      \gamma'\right>.
  \end{gather}

  For $\gamma, \gamma' \in \hBB'$, we have the following (combinatorial)
  formula:
  \begin{equation}
    \label{equ:scpandpartial}
    \left<\exp(\partial)(\gamma\gamma'), 1\right>
    = \left<\exp(\partial)\gamma, \exp(\partial)\gamma'\right>.
  \end{equation}
\end{prop}

\begin{proof}
  The formula~\eqref{equ:ellandpartial} should be clear from the
  definitions.

  Let us investigate~\eqref{equ:scpandpartial} a bit more.
  We can assume that $\gamma$ and $\gamma'$ are Jacobi diagrams with
  $l$ resp.\ $l'$ univalent vertices and $l + l' = 2n$ with $n \in \N_0$. So
  we have to prove
  \begin{equation*}
    \label{equ:combin}
    \frac{\partial^n}{n!} (\gamma\gamma')
    = \sum_{\substack{m, m' = 0 \\ l - 2 m = l' - 2 m'}}^\infty
      \left<\frac{\partial^m}{m!} \gamma, \frac{\partial^{m'}}{m'!}
        \gamma'\right>,
  \end{equation*}
  since $\left<\cdot, 1\right>: \hBB \to \hBBo$ means to remove the components
  with at least one univalent vertex. Recalling the meaning of $\left<\cdot,
    \cdot\right>$, it should be clear that~\eqref{equ:scpandpartial} follows
  from the fact that applying $\frac{\partial^k}{k!}$ on a Jacobi diagram
  means to glue all subsets of $2k$ of its univalent vertices to $k$ pairs in
  all possible ways.
\end{proof}

\subsection{An $\sltwo$-action on the space of graph homology}
\label{sec:sltwo}

In this short section we want to extend the space of graph homology slightly.
This is mainly due to two reasons: When we defined the expression
$\hat\Gamma(\Gamma)$ for two Jacobi diagrams $\Gamma$ and $\Gamma'$, we
restricted ourselves to the case that $\Gamma$ or $\Gamma'$ does not contain a
component with an $\ell$. Secondly, we have not given the \emph{zero-wheel
  $w_0$} a meaning yet.

We do this by adding an element $\bigcirc$ to the various spaces of graph
homology.

\begin{defn}
  \emph{The extended space of graph homology} is the space
  $\hBB[[\bigcirc]]$. Further, we set $w_0 := \bigcirc$, which, at least
  picturally, is in accordance with the definition of $w_k$ for $k > 0$.
\end{defn}

Note that this element is not depicting a Jacobi diagram as we have defined
it. Nevertheless, we want to use the notion that $\bigcirc$ has no univalent
and no trivalent vertices, i.e.\ the homogeneous component of degree zero of
$\hBB[[\bigcirc]]$ is $\C[[\bigcirc]]$.

When defining $\Gamma/(u, u')$ for a Jacobi diagramm $\Gamma$ with two
univalent vertices $u$ and $u'$, i.e.\ gluing $u$ to $u'$, we assumed that $u$
and $u'$ are not the vertices of one component $\ell$ of $\Gamma$. Now we
extend this definition by defining $\Gamma/(u, u')$ to be the extended graph
homology class we get by replacing $\ell$ with $\bigcirc$, whenever $u$ and
$u'$ are the two univalent vertices of a component $\ell$ of $\Gamma$.

Doing so, we can give the expression $\hat\gamma(\gamma') \in
\hBB[[\bigcirc]]$ a meaning with no restrictions on the two graph homology
classes $\gamma, \gamma' \in \hBB$, i.e.\ every $\gamma \in \hBB[[\bigcirc]]$
defines a $\hBBo[[\bigcirc]]$-linear map
\begin{equation}
  \hat\gamma: \hBB[[\bigcirc]] \to \hBB[[\bigcirc]].
\end{equation}

\begin{exmp}
  We have
  \begin{equation}
    \partial\ell = \bigcirc.
  \end{equation}
\end{exmp}

\begin{rmk}
  We can similarly extend $\left<\cdot, \cdot\right>: \hBB' \times \hBB \to
  \hBBo$ to a $\hBBo[[\bigcirc]]$-bilinear form
  \begin{equation}
    \left<\cdot, \cdot\right>: \hBB[[\bigcirc]] \times \hBB[[\bigcirc]] \to
    \hBBo[[\bigcirc]].
  \end{equation}
\end{rmk}

Both $\ell/2$ and $\partial$ are two operators acting on the extended space of
graph homology, the first one just multiplication with $\ell/2$. By
calculating their commutator, we show that they induce a natural structure of
an $\sltwo$-module on $\hBB[[\bigcirc]]$.

\begin{prop}
  Let $H: \hBB[[\bigcirc]] \to \hBB[[\bigcirc]]$ be the linear operator which
  acts on $\gamma \in \hBB_{k, l}[[\bigcirc]]$ by
  \begin{equation}
    H \gamma = \left(\frac 1 2 \bigcirc + l\right) \gamma.
  \end{equation}

  We have the following commutator relations in $\End\hBB[[\bigcirc]]$:
  \begin{align}
    [\ell/2, \partial] & = - H,
    \label{equ:sltwocomm}
    \\
    [H, \ell/2] & = 2 \cdot \ell/2,
    \label{equ:hell}
    \\
    \intertext{and}
    [H, \partial] & = - 2 \partial,
    \label{equ:hpartial}
  \end{align}
  i.e.\ the triple $(\ell/2, - \partial, H)$ defines a $\sltwo$-operation on
  $\hBB[[\bigcirc]]$.
\end{prop}

\begin{proof}
  Equations~\eqref{equ:hell} and~\eqref{equ:hpartial} follow from the fact
  that multiplying by $\bigcirc$ commutes with $\ell/2$ and $\partial$, and
  from the fact that $\ell/2$ is an operator of degree $2$ with respect to the
  grading given by the number of univalent vertices, whereas $\partial$ is an
  operator of degree $-2$ with respect to the same grading.

  It remains to look at~\eqref{equ:sltwocomm}. For $\gamma \in \hBB_{k,
    l}[[\bigcirc]]$, we calculate
  \begin{gather}
    [\ell, \partial]\gamma = \ell \partial(\gamma) - \partial (\ell \gamma)
    = \ell \partial(\gamma) - \partial(\ell) \gamma - \ell \partial(\gamma)
    - \partial(\ell, \gamma)
    = - \bigcirc \gamma - 2 l \gamma
    = - 2 H \gamma.
  \end{gather}
\end{proof}

\begin{rmk}
  Since $\hBB[[\bigcirc]]$ is infinite-dimensional, we have unfortunately
  difficulties to apply the standard theory of $\sltwo$-representations to
  this $\sltwo$-module. For example, there are no eigenvectors for the
  operator $H$.
\end{rmk}

\subsection{Closed and connected graphs, the closure of a graph}
As the number of connected components of a Jacobi diagram is preserved by the
IHX- and AS-relations each graph homology space inherits a grading by the
number of connected components. For any $k \in \N_0$ we define $\BB^k$ to be
the subspace of $\BB$ spanned by all Jacobi diagrams with exactly $k$
connected components. Similarly, we define $\BBo^k$, $\hBB^k$, $\hBBo^k$.

We have $\BB = \bigoplus_{k = 0}^\infty \BB^k$ with $\BB^0 = \Q \cdot
1$. Analogous results hold for $\BBo$, $\hBB$, $\hBBo$.

\begin{defn}
  A graph homology class $\gamma$ is called \emph{closed} if $\gamma \in
  \hBBo$. The class $\gamma$ is called \emph{connected} if $\gamma \in
  \hBB^1$. The \emph{connected component of $\gamma$} is defined to be
  $\pr^1(\gamma)$ where $\pr^1: \hBB = \prod_{i = 0}^\infty \hBB^i \to \hBB^1$
  is the canonical projection. The \emph{closure $\clos\gamma$ of $\gamma$} is
  defined by $\clos\gamma := \clos{\gamma, \exp(\ell/2)}$.  The
  \emph{connected closure $\cclos\gamma$ of $\gamma$} is defined to be the
  connected component of the closure $\clos\gamma$ of $\gamma$.
\end{defn}

For every finite set $L$, we define $\partitions_2(L)$ to be the set of
partitions of $L$ into subsets of two elements. With this definition, we can
express the closure of a Jacobi diagram $\Gamma$ as
\begin{gather}
  \clos{\Gamma} = \sum_{\pi \in \partitions_2(L)} \Gamma/\pi.
\end{gather}

\begin{exmp}
  We have $\clos{w_2} = \ThetaGraph$, $\cclos{w_2} = \ThetaGraph$,
  $\clos{w_2^2} = 2 \ThetaGraph_2 + \ThetaGraph^2$, $\cclos{w_2^2} = 2
  \ThetaGraph_2$.
\end{exmp}

Let $L_1, \dots, L_n$ be finite and pairwise disjoint sets. We set $L :=
\bigsqcup_{i = 1}^n L_i$. Let $\pi \in \partitions_2(L)$ be a partition of $L$
in 2-element-subsets. We say that a pair $l, l' \in L$ \emph{is linked by
  $\pi$} if there is an $i \in \set{1, \dots n}$ such that ${l, l'} \in L_i$
or $\set{l, l'} \in \pi$. We say that \emph{$\pi$
  connects the sets $L_1, \dots, L_n$} if and only if for each pair $l, l'
\in L$ there is a chain of elements $l_1, \dots, l_k$ such that $l$ is linked
to $l_1$, $l_i$ is linked to $l_{i + 1}$ for $i \in \set{1, \dots, k - 1}$ and
$l_k$ is linked to $l'$. The subset of $\partitions_2(L)$ of partitions $\pi$
connecting $L_1, \dots, L_n$ is denoted by $\partitions_2(\set{L_1, \dots,
  L_n})$.
We have
\begin{gather}
  \label{eq:partdec}
  \partitions_2(L) = \bigsqcup_{\bigsqcup \mathfrak I = \set{1, \dots, n}}
  \set{\bigsqcup_{I \in \mathfrak I} \pi_I : \pi_I \in P_2(\set{L_i: i \in
  I})}.
\end{gather}
Here, $\bigsqcup \mathfrak I = \set{1, \dots, n}$ means that $\mathfrak I$
is a partition of $\set{1, \dots, n}$ in disjoint subsets.

Let $\Gamma_1, \dots, \Gamma_n$ be connected Jacobi diagrams. We denote by
$\Gamma := \prod_{i = 1}^n \Gamma_i$ the product over all these Jacobi
diagrams. Let $L_i$ be the set of legs of $\Gamma_i$ and denote by $L :=
\bigsqcup_{i = 1}^n L_i$ the set of all legs of $\Gamma$.

For every partition $\pi \in \partitions_2(L)$ the graph
$\Gamma/\pi$ is connected if and only if $\pi \in \partitions_2(\set{L_1,
  \dots, L_n})$.

Using~\eqref{eq:partdec} we have
\begin{multline}
  \label{eq:closofprod}
  \clos\Gamma = \sum_{\pi \in \partitions_2(L)} \Gamma/\pi
  = \sum_{\bigsqcup \mathfrak I = \set{1, \dots, n}} \prod_{I \in \mathfrak I}
  \sum_{\pi \in \partitions_2(\set{L_i: i \in I})}
  \left(\prod_{i \in I} \Gamma_i\right)/\pi
  \\
  = \sum_{\bigsqcup \mathfrak I = \set{1, \dots, n}} \prod_{I \in \mathfrak I}
  \cclos{\prod_{i \in I} \Gamma_i}.
\end{multline}

With this result we can prove the following Proposition:
\begin{prop}
  \label{prop:cclos}
  For any connected graph homology class $\gamma$ we have
  \begin{gather}
    \label{eq:cclos}
    \exp\cclos{\exp\gamma} = \clos{\exp{\gamma}}.
  \end{gather}
\end{prop}

Note that both sides are well-defined in $\hBB$ since $\gamma$ and
$\cclos{\cdots}$ as connected graphs have no component in degree zero.

\begin{proof}
  Let $\Gamma$ be any connected Jacobi diagram. By~\eqref{eq:closofprod} we
  have
  \begin{gather*}
    \clos{\Gamma^n} = \sum_{\bigsqcup \mathfrak I = \set{1, \dots, n}}
    \prod_{I \in \mathfrak I} \cclos{\Gamma^{\# I}}
    = \sum_{\lambda \in \partitions(n)} n! \prod_{i = 1}^\infty
    \frac 1 {\lambda_i!} \left(\cclos{\Gamma^i}/i!\right)^{\lambda_i}.
  \end{gather*}
  By linearity this result holds also if we substitute $\Gamma$ by the
  connected graph homology class $\gamma$.
  
  Using this,
  \begin{multline*}
    \clos{\exp\gamma} = \sum_{n = 0}^\infty \frac 1 {n!} \clos{\gamma^n}
    = \sum_{n = 0}^\infty \frac 1 {n!}
    \sum_{\lambda \in \partitions(n)} n! \prod_{i = 1}^\infty
    \frac 1 {\lambda_i!} \left(\cclos{\gamma^i}/i!\right)^{\lambda_i}
    \\
    = \prod_{i = 1}^\infty \sum_{\lambda = 0}^\infty \frac 1{\lambda!}
    \left(\cclos{\gamma^i}/i!\right)^\lambda
    = \prod_{i = 1}^\infty \exp\left(\cclos{\gamma^i}/i!\right)
    = \exp\cclos{\exp\gamma}.
  \end{multline*}
\end{proof}

\subsection{Polywheels}

\begin{defn}
  For each $n \in \N_0$ we set $\nwheel_{2n} := - w_{2n}$. Let $\lambda$ be a
  partition of $n$. We set
  \begin{gather}
    \nwheel_{2 \lambda} := \prod_{i = 1}^\infty \nwheel_{2i}^{\lambda_i}.
  \end{gather}
  The closure $\clos{\nwheel_{2 \lambda}}$ of $\nwheel_{2 \lambda}$ is called
  a \emph{polywheel}. The subspace in $\BBo$ spanned by all polywheels is
  denoted by $\mathcal W$ and called the \emph{polywheel subspace}. The
  subalgebra in $\BBo$ spanned by all polywheels is denoted by $\mathcal C$
  and called the \emph{algebra of polywheels}.
  
  The connected closure $\cclos{\nwheel_{2 \lambda}}$ of $\nwheel_{2 \lambda}$
  is called a \emph{connected polywheel}.
\end{defn}

\begin{rmk}
  As discussed by J.~Sawon in his thesis~\cite{sawon99}, $\mathcal W$ is
  proper graded subspace of $\mathcal \BBo$. From degree eight on, $\BBo_k$
  is considerably larger than $\mathcal W_k$. On the other hand it is unknown
  (at least to the author) if the inclusion $\mathcal C \subseteq \BBo$ is
  proper.
\end{rmk}

\begin{rmk}
 The subalgebra $\mathcal C'$ in $\BBo$ spanned by all connected polywheels
 equals $\mathcal C$. This is since we can use~\eqref{eq:cclos} to express
 every polywheel as a polynomial of connected polywheels and vice versa.
\end{rmk}

\begin{exmp}
  Using Proposition~\ref{prop:cclos} we calculated the following expansions of
  the connected polywheels in terms of wheels:
  \begin{gather}
    \begin{aligned}
      \cclos{\nwheel_2} & = \clos{\nwheel_2}
      \\[1ex]
      \cclos{\nwheel_2^2} & = \clos{\nwheel_2^2} - \clos{\nwheel_2}^2 \\
      \cclos{\nwheel_4} & = \clos{\nwheel_4}
      \\[1ex]
      \cclos{\nwheel_2^3} & = \clos{\nwheel_2^3} - 3 \clos{\nwheel_2}
      \clos{\nwheel_2^2} + 2 \clos{\nwheel_2}^3 \\
      \cclos{\nwheel_2 \nwheel_4} & = \clos{\nwheel_2 \nwheel_4}
      - \clos{\nwheel_2} \clos{\nwheel_4} \\
      \cclos{\nwheel_6} & = \clos{\nwheel_6}
      \\[1ex]
      \cclos{\nwheel_2^4} & = \clos{\nwheel_2^4}
      - 4 \clos{\nwheel_2} \clos{\nwheel_2^3} - 3 \clos{\nwheel_2^2}^2
      + 12 \clos{\nwheel_2}^2 - 6 \clos{\nwheel_2}^4 \\
      \cclos{\nwheel_2^2 \nwheel_4} & = \clos{\nwheel_2^2 \nwheel_4}
      - 2 \clos{\nwheel_2} \clos{\nwheel_2 \nwheel_4}
      - \clos{\nwheel_2^2} \clos{\nwheel_4}
      + 2 \clos{\nwheel_2}^2 \clos{\nwheel_4} \\
      \cclos{\nwheel_2 \nwheel_6} & = \clos{\nwheel_2 \nwheel_6}
      - \clos{\nwheel_2} \clos{\nwheel_6} \\
      \cclos{\nwheel_4^2} & = \clos{\nwheel_4^2} - \clos{\nwheel_4}^2 \\
      \cclos{\nwheel_8} &= \clos{\nwheel_8}.
    \end{aligned}
  \end{gather}
\end{exmp}

\section{Holomorphic symplectic manifolds}
\label{sec:sympl}

\subsection{Definition and general properties}

\begin{defn}
  A \emph{holomorphic symplectic manifold $(X, \sigma)$} is a compact complex
  manifold $X$ together with an everywhere non-degenerate holomorphic two-form
  $\sigma \in \HH^0(X, \Omega_X^2)$. Here, we call $\sigma$ \emph{everywhere
    non-degenerate} if $\sigma$ induces an isomorphism $\Tang_X \to \Omega_X$.

  The holomorphic symplectic manifold $(X, \sigma)$ is called
  \emph{irreducible} if it is simply-connected and $\HH^0(X, \Omega_X^2)$ is
  one-dimensional, i.e.~spanned by $\sigma$.
\end{defn}

It follows immediately that every holomorphic symplectic manifold $X$ has
trivial canonical bundle whose sections are multiples of $\sigma^n$, and,
therefore, vanishing first Chern class. In fact, all odd Chern classes
vanish:

\begin{prop}
  Let $X$ be a complex manifold and $E$ a complex vector bundle on $X$. If $E$
  admits a symplectic two-form, i.e.~there exists a section $\sigma \in
  \HH^0(X, \Lambda^2 E^*)$ such that the induced morphism $E \mapsto E^*$ is
  an isomorphism, all odd Chern classes of $E$ vanish.
\end{prop}

\begin{rmk}
  That the odd Chern classes of $E$ vanish up to two-torsion follows
  immediately from the fact $c_{2k + 1}(E) = - c_{2k + 1}(E^*)$ for $k \in
  \N_0$.
\end{rmk}

The following proof using the splitting principle has been suggested to me by
Manfred Lehn.
\begin{proof}
  We prove the proposition by induction over the rank of $E$. For $\rk E = 0$,
  the claim is obvious.
  
  By the splitting principle (see e.g.~\cite{fulton96}), we can assume that
  $E$ has a subbundle $L$ of rank one. Let $L^\orth$ be the
  $\sigma$-orthogonal subbundle to $L$ of $E$. Since $\sigma$ is symplectic,
  $L^\orth$ is of rank $n-1$ and $L$ is a subbundle of $L^\orth$. We have the
  following short exact sequences of bundles on $X$:
  \begin{gather*}
    \begin{CD}
      0 @>>> L @>>> E @>>> E/L @>>> 0
    \end{CD}
    \intertext{and}
    \begin{CD}
      0 @>>> L^\orth/L @>>> E/L @>>> E/L^\orth @>>> 0.
    \end{CD}
  \end{gather*}
  Since $\sigma$ induces a symplectic form on $L^\orth/L$, by induction, all
  odd Chern classes of this bundle of rank $\rk E - 2$ vanish. Furthermore,
  note that $\sigma$ induces an isomorphism between $L$ and $(E/L^\orth)^*$, so
  all odd Chern classes of $L \oplus E/L^\orth$ vanish.

  Now, the two exact sequences give us $c(E) = c(L \oplus E/L^\orth) \cdot
  c(L^\orth/L)$. Therefore, we can conclude that all odd Chern classes of $E$
  vanish.
\end{proof}

\begin{prop}
  For any irreducible holomorphic symplectic manifold $(X, \sigma)$ of
  dimension $2n$ and $k \in {0, \dots, n}$ the space $\HH^{2k}(X, \OX)$
  is one-dimensional and spanned by the cohomology class $[\bar\sigma]^k$.
\end{prop}

\begin{proof}
  See~\cite{beauville83}.
\end{proof}

\subsection{A pairing on the cohomology of a holomorphic symplectic manifold}

Let $(X, \sigma)$ be a holomorphic symplectic manifold.
There is a natural pairing of coherent sheafs
\begin{gather}
  \Lambda_* \Tang_X \otimes \Lambda^* \Omega_X \to \OX.
\end{gather}
As the natural morphism from $\Lambda_* \Tang_X$ to $\Lambda^* \Tang_X$ is an
isomorphism and $\Lambda^* \Tang_X$ can be identified with $\Lambda^*
\Omega_X$ by means of the symplectic form, we therefore have a natural map 
\begin{gather}
  \Lambda^* \Omega \otimes \Lambda^* \Omega_X \to \OX.
\end{gather}
We write
\begin{gather}
  \scal{\cdot, \cdot}: \HH^p(X, \Omega^*) \otimes \HH^q(X, \Omega^*) \to
  \HH^{p + q}(X, \OX), (\alpha, \beta) \mapsto \scal{\alpha, \beta}
\end{gather}
for the induced map for any $p, q \in \N_0$.

In~\cite{nieper02} we proved the following proposition:

\begin{prop}
  For any $\alpha \in \HH^*(X, \Omega^*)$ we have
  \begin{gather}
    \int_X \alpha \exp\sigma = \int_X \scal{\alpha, \exp \sigma} \exp
    \sigma.
  \end{gather}
\end{prop}

\subsection{Example series}

There are two main series of examples of irreducible holomorphic symplectic
manifolds. Both of them are based on the Hilbert schemes of points on a
surface:

Let $X$ be any smooth projective surface over $\C$ and $n \in
\N_0$. By $\Hilb n
X$ we denote the Hilbert scheme of zero-dimensional subschemes of length $n$
of $X$. By a result of Fogarty~(\cite{fogarty68}), $\Hilb n X$ is a smooth
projective variety of dimension $2n$. The Hilbert scheme can be viewed as a
resolution $\rho: \Hilb n X \to X^{(n)}$ of the $n$-fold symmetric product
$X^{(n)} := X^n/\mathfrak S_n$. The morphism $\rho$, sending closed points,
i.e.~subspaces of $X$, to their support counting multiplicities, is called the
Hilbert-Chow morphism.

Let $\alpha \in \HH^2(X, \C)$ be any class. The class $\sum_{i = 1}^n \pr_i^*
\alpha \in \HH^2(X^n, \C)$ is invariant under the action of $\mathfrak S_n$,
where $\pr_i: X^n \to X$ denotes the projection on the $i^{\mathrm{th}}$
factor. Therefore, there exists a class $\alpha^{(n)} \in \HH^2(X^{(n)}, \C)$
with $\pi^* \alpha^{(n)} = \sum_{i = 1}^n \pr_i^* \alpha$, where $\pi: X^n \to
X^{(n)}$ is the canonical projection. Using $\rho$ this induces a class $\Hilb
n \alpha$ in $\HH^2(\Hilb n X, \C)$.

If $X$ is a K3 surface or an abelian surface, there exists a holomorphic
symplectic form $\sigma \in \HH^{2, 0}(X) \subseteq \HH^2(X, \C)$. It was shown
by Beauville in~\cite{beauville83} that $\Hilb n \sigma$ is again symplectic,
so $(\Hilb n X, \Hilb n \sigma)$ is a holomorphic symplectic
manifold.

\begin{exmp}
  For any K3 surface $X$ and holomorphic symplectic form $\sigma \in \HH^{2,
  0}(X)$, the pair $(\Hilb n X, \Hilb n \sigma)$ is
  in fact an irreducible holomorphic symplectic manifold.
\end{exmp}
This has also been proven by Beauville.  In the case of an abelian surface
$A$, we have to work a little bit more as $\Hilb n A$ is not irreducible in
this case:

Let $A$ be an abelian surface and let us denote by $s: \Hilb n A \to A$ the
composition of the summation morphism $A^{(n)} \to A$ with the Hilbert-Chow
morphism $\rho: \Hilb n A \to A^{(n)}$.

\begin{defn}
  For any $n \in \N$, the \emph{$n^{\mathrm{th}}$ generalised Kummer variety
    $\Kummer n A$} is the fibre of $s$ over $0 \in A$.
  For any class $\alpha \in \HH^2(A, \C)$, we set $\Kummer n \alpha := \Hilb
  n \alpha|_{\Kummer n A}$.
\end{defn}

\begin{rmk}
  For $n = 2$ the generalised Kummer variety coincides with the Kummer model
  of a K3 surface (therefore the name).
\end{rmk}

\begin{exmp}
  For every abelian surface $A$ and holomorphic symplectic form $\sigma \in
  \HH^{2, 0}(A)$, the pair $(\Kummer n A, \Kummer n \sigma)$ is an irreducible
  holomorphic symplectic manifold of dimension $2n - 2$.
\end{exmp}
The proof can also be found in~\cite{beauville83}.

\subsection{About $\Hilb n \alpha$ and $\Kummer n \alpha$}

Let $X$ be any smooth projective surface and $n \in \N_0$.

Let $\Hilb {n, n + 1} X$ denote the incidence variety of all pairs $(\xi,
\xi') \in \Hilb n X \times \Hilb {n + 1} X$ with $\xi \subseteq \xi'$
(see~\cite{lehn01}). We denote by $\psi: \Hilb {n, n + 1} X \to \Hilb {n + 1}
X$ and by $\phi: \Hilb {n, n + 1} X \to \Hilb n X$ the canonical maps. There
is a third canonical map $\chi: \Hilb {n, n + 1} X \to X$ mapping $(\xi, \xi')
\mapsto x$ if $\xi'$ is obtained by extending $\xi$ at the closed point $x \in
X$.

\begin{prop}
  \label{prop:psialpha}
  For any $\alpha \in \HH^2(X, \C)$ we have
  \begin{gather}
    \psi^* \Hilb{n + 1} \alpha = \phi^* \Hilb n \alpha + \chi^* \alpha.
  \end{gather}
\end{prop}
 
\begin{proof}
  Let $p: X^{(n)} \times X \to X^{(n)}$ and $q: X^{(n)} \times X \to X$ denote
  the canonical projections. Let $\tau: X^{(n)} \times X \to X^{(n + 1)}$ the
  obvious symmetrising map. 
  The following diagram
  \begin{gather*}
    \begin{CD}
      \Hilb {n, n + 1} X @= \Hilb {n, n + 1} X \\
      @V{(\phi, \chi)}VV @VV{\psi}V \\
      \Hilb n X \times X & & \Hilb {n + 1} X \\
      @V{\rho \times \id_X}VV @VV{\rho'}V \\
      X^{(n)} \times X @>{\tau}>> X^{(n + 1)} \\
      @A{\pi \times \id_X}AA @AA{\pi'}A \\
      X^{n + 1} @= X^{n + 1}
    \end{CD}
  \end{gather*}
  is commutative. (Note that we have primed some maps to avoid name
  clashes.) We claim that $\tau^* \alpha^{(n + 1)} = p^*
  \alpha^{(n)} + q^* \alpha$. In fact, since
  \begin{gather*}
    (\pi \times \id_X)^*\tau^*
    \alpha^{(n + 1)} = \pi^{\prime *} \alpha^{(n + 1)}
    = \sum_{i = 1}^{n + 1} \pr_i^* \alpha,
  \end{gather*}
  this follows from the definition of $\alpha^{(n)}$.
  Finally, we can read off the diagram that
  \begin{multline*}
    \psi^* \Hilb{n + 1} \alpha = \psi^* \rho^{\prime *} \alpha^{(n + 1)}
    = (\phi, \chi)^* (\rho \times \id_X)^* \tau^* \alpha^{(n + 1)}
    \\
    = (\phi, \chi)^* (\rho \times \id_X)^* (p^* \alpha^{(n)} + q^* \alpha)
    = \phi^* \Hilb n \alpha + \chi^* \alpha.
  \end{multline*}
\end{proof}

\begin{prop}
  \label{prop:splitalpha}
  Let $X = X_1 \sqcup X_2$ be the disjoint union of two projective smooth
  surfaces $X_1$ and $X_2$. We then have
  \begin{gather}
    \Hilb n X = \bigsqcup_{n_1 + n_2 = n} \Hilb {n_1} X_1 \times \Hilb {n_2}
    X_2.
  \end{gather}
  If $\alpha \in \HH^2(X, \C)$ decomposes as $\alpha|_{X_1} = \alpha_1$ and
  $\alpha|_{X_2} = \alpha_2$, then $\Hilb n \alpha$ decomposes as
  \begin{gather}
    \Hilb n \alpha|_{\Hilb {n_1} X_! \times \Hilb {n_2} X_2}
    = \pr_1^* \Hilb {n_1} \alpha_1 + \pr_2^* \Hilb {n_2} \alpha_2.
  \end{gather}
\end{prop}

\begin{proof}
  The splitting of $\Hilb n X$ follows from the universal property of the
  Hilbert scheme and is a well-known fact. The statement on $\Hilb n \alpha$
  is easy to prove and so we shall only give a sketch: Let us denote by $i:
  \Hilb{n_1} X_1 \times \Hilb{n_2} X_2 \to \Hilb n X$ the natural inclusion.
  Furthermore let $j: X_1^{(n_1)} \times X_2^{(n_2)} \to X^{(n)}$ denote the
  natural symmetrising map. The following diagram is commutative:
  \begin{gather}
    \begin{CD}
      \Hilb{n_1} X_1 \times \Hilb{n_2} X_2 @>i>> \Hilb n X \\
      @V{\rho_1 \times \rho_2}VV @VV{\rho}V \\
      X_1^{(n_1)} \times X_2^{(n_2)} @>>{j}> X^{(n)},
    \end{CD}
  \end{gather}
  where the $\rho_i: \Hilb{n_i} X_i \to X_i^{(n_i)}$ are the
  Hilbert-Chow morphisms. Since $j^* \alpha^{(n)} = \pr_1^* \alpha_1^{(n_1)}
  + \pr_2^* \alpha_2^{(n_2)}$, the commutativity of the diagram proves the
  statement on $\Hilb n \alpha$.
\end{proof}

Let $A$ be again an abelian surface and $n \in \N$. Since $A$ acts on
itself by translation, there is also an induced operation of $A$ on the
Hilbert scheme $\Hilb n A$. Let us denote the restriction of this operation to
the generalised Kummer variety $\Kummer n A$ by $\nu: A \times \Kummer n A \to
\Hilb n A$. It fits into the following cartesian square:
\begin{gather}
  \begin{CD}
    A \times \Kummer n A @>{\nu}>> \Hilb n A \\
    @V{\pr_1}VV @VV{s}V \\
    A @>>{n}> A,
  \end{CD}
\end{gather}
where $s$ is the summation map as having been defined above and $n: A \to A, a
\mapsto n a$ is the (multiplication-by-$n$)-morphism. Since $n$ is a Galois
cover of degree $n^4$, the same holds true for $\nu$.

\begin{prop}
  \label{prop:pullbackofalpha}
  For any $\alpha \in \HH^2(A, \C)$, we have
  \begin{gather}
    \nu^* \Hilb n \alpha = n \pr_1^* \alpha + \pr_2^* \Kummer n \alpha.
  \end{gather}
\end{prop}

\begin{proof}
  By the K\"unneth decomposition theorem, we know that $\nu^* \Hilb n \alpha$
  splits:
  \begin{gather*}
    \nu^* \Hilb n \alpha = \pr_1^* \alpha_1 + \pr_2^* \alpha_2.
  \end{gather*}
  Set $\iota_1: A \to A \times
  \Kummer n A, a \mapsto (a, \xi_0)$ and $\iota_2: \Kummer n A \to A \times
  \Kummer n A, \xi \mapsto (0, \xi)$, where $\xi_0$ is any subscheme of length
  $n$ concentrated in $0$. We have
  \begin{gather}
    \alpha_1 = \iota_1^*\nu^*\Hilb n \alpha
    = (\rho \circ \nu \circ \iota_1)^* \alpha^{(n)}
    = (a \mapsto (\underbrace{a, \dots a}_n))^* \alpha^{(n)} = n \alpha
  \end{gather}
  and
  \begin{gather}
    \alpha_2 = \iota_2^*\nu^*\Hilb n \alpha = i^* \Hilb n \alpha = \Kummer n
    \alpha, 
  \end{gather}
  where $i: \Kummer n A \to \Hilb n A$ is the natural inclusion map, thus
  proving the proposition.
\end{proof}

\subsection{Complex genera of Hilbert schemes of points on surfaces}

The following theorem is an adaption of Theorem~4.1 of~\cite{lehn01} to our
context.
\begin{thm}
  \label{thm:ptoptilde}
  Let $P$ be a polynomial in the variables $c_1, c_2, \dots$ and $\alpha$ over
  $\Q$. There exists a polynomial $\tilde P \in \Q[z_1, z_2, z_3, z_4]$ such
  that for every smooth projective surface $X$, $\alpha \in \HH^2(X, \Q)$ and
  $n \in \N_0$ we have:
  \begin{gather}
    \int_{\Hilb n X} P(c_*(\Hilb n X), \Hilb n \alpha)
    = \tilde P \left(\int_X \alpha^2/2, \int_X c_1(X) \alpha, \int_X
      c_1(X)^2/2, \int_X c_2(X)\right).
  \end{gather}
\end{thm}

\begin{proof}
  The prove goes along the very same lines as the proof of Proposition~0.5
  in~\cite{lehn01} (see there). The only new thing we need is
  Proposition~\ref{prop:psialpha} of this paper to be used in the induction
  step of the adapted proof of Proposition~3.1 of~\cite{lehn01} to our
  situation.
\end{proof}

Let $R$ be any $\Q$-algebra (commutative and with unit) and
let $\phi \in R[[c_1, c_2, \dots]]$ be a non-vanishing power series in the
universal Chern classes such that $\phi$ is multiplicative with respect to the
Whitney sum of vector bundles, i.e.
\begin{gather}
  \phi(E \oplus F) = \phi(E) \phi(F)
\end{gather}
for all complex manifolds and complex vector bundles $E$ and $F$ on $X$. Any
$\phi$ with this property induces a complex genus, also denoted by $\phi$, by
setting $\phi(X) := \int_X \phi(\Tang_X)$ for $X$ a compact complex
manifold. Let us call such a $\phi$ \emph{multiplicative}.

\begin{rmk}
  By Hirzebruch's theory of multiplicative sequences and complex genera
  (\cite{hirzebruch56}), we know that
  \begin{enumerate}
  \item each complex genus is induced by a unique multiplicative $\phi$, and
  \item the multiplicative elements in $R[[c_1, c_2, \dots]]$ are exactly
    those of the form $\exp(\sum_{k = 1}^\infty a_k s_k)$ with $a_k \in R$.
  \end{enumerate}
\end{rmk}

More or less formally the following theorem follows from
Theorem~\ref{thm:ptoptilde}.
\begin{thm}
  \label{thm:abcdphi}
  For each multiplicative $\phi \in R[[c_1, c_2, \dots]]$, there exist unique
  power series $A_\phi(p), B_\phi(p), C_\phi(p), D_\phi(p) \in pR[[p]]$ with
  vanishing constant coefficient such that for all smooth projective surfaces
  $X$ and $\alpha \in \HH^2(X, \C)$ we have:
  \begin{multline}
    \label{eq:abcdphi}
    \sum_{n = 0}^\infty \left(\int_{\Hilb n X} \phi(\Hilb n X)\exp(\Hilb n
      \alpha) \right) p^n \\
    = \exp\left( A_\phi(p) \int_X \alpha^2/2 +
      B_\phi(p) \int_X c_1(X) \alpha\right.
    \\
    + \left.C_\phi(p) \int_X c_1^2(X)/2
      + D_\phi(p) \int_X c_2(X) \right).
  \end{multline}
  The first terms of $A_\phi(p), B_\phi(p), C_\phi(p), D_\phi(p)$
  are given by
  \begin{multline}
    A_\phi(p) = p + O(p^2), \quad
    B_\phi(p) = \phi_1 p + O(p^2), \quad
    C_\phi(p) = \phi_{11} p + O(p^2), \text{and}\\
    D_\phi(p) = \phi_2 p + O(p^2),
  \end{multline}
  where $\phi_1$ is the coefficient of $c_1$ in $\phi$, $\phi_{11}$ the
  coefficient of $c_1^2/2$ and $\phi_2$ the coefficient of $c_2$.
\end{thm}

\begin{proof}
  This theorem is again an adaption of a theorem~(Theorem~4.2)
  of~\cite{lehn01} to our context. Nevertheless, let us give the proof here:
  
  Set $K := \set{(X, \alpha): \text{$X$ is a smooth projective surface and
      $\alpha \in \HH^2(X, \C)$}}$ and let $\gamma: K \to \Q^4$ be the map
  $(X, \alpha) \mapsto (\alpha^2/2, c_1(X) \alpha, c_1(X)^2/2, c_2(X))$. Here,
  we have supressed the integral signs $\int_X$ and interpret the expressions
  $\alpha^2$, etc.\ as intersection numbers on $X$. The image of $K$ spans the
  whole $\Q^4$ (for explicit generators, we refer to~\cite{lehn01}).
  
  Now let us assume that a $(X, \alpha) \in K$ decomposes as $(X, \alpha)
  = (X_1, \alpha_1) \sqcup (X_2, \alpha_2)$. By the multiplicative behaviour
  of $\phi$ and $\exp$ we see that
  \begin{multline*}
     \int_{\Hilb n X} \phi(c_*(\Hilb n X)) \exp(\Hilb n \alpha)
     \\
     = \sum_{n_1 + n_2 = n}
     \left(
       \int_{\Hilb {n_1} X_1} \phi(c_*(\Hilb {n_1} X)) \exp(\Hilb {n_1}
       \alpha_2)\right)
     \left(
       \int_{\Hilb {n_2} X_2} \phi(c_*(\Hilb {n_2} X)) \exp(\Hilb {n_2}
       \alpha_2)\right),
  \end{multline*}
  whereas $H_\phi(p)(X, \alpha) := 
  \sum_{n = 0}^\infty \left(\int_{\Hilb n X} \phi(\Hilb n X)\exp(\Hilb n
    \alpha) \right) p^n$ fulfills
  \begin{gather*}
    \tag{$*$}
    H_\phi(p)(X, \alpha) = H_\phi(p)(X_1, \alpha_1) H_\phi(p)(X_2, \alpha_2).
  \end{gather*}
  Since $H_\phi(p): K \to \Q^4$ factors through $\gamma$ and a map $h: \Q^4
  \to R[[p]]$ by Theorem~\ref{thm:ptoptilde} and as the image of $\gamma$ is
  Zariski dense in $\Q^4$, we conclude from $(*)$ that $\log h$ is a linear
  function which proves the first part of the theorem.

  To get the first terms of the power series, we expand both sides
  of~\eqref{eq:abcdphi}. The left hand side expands as
  \begin{gather}
    1 + (\alpha^2/2 + \phi_1 c_1(X) \alpha + \frac{\phi_{11}} 2 c_1^2(X) +
    \phi_2 c_2(X)) p + O(p^2),
  \end{gather}
  while the right hand side expands as
  \begin{gather}
    1 + (A_1 \alpha^2/2 + B_1 c_1(X) \alpha + C_1 c_1(X)^2 + D_1 c_2(X)) p
    + O(p^2),
  \end{gather}
  where $A_1, B_1, C_1, D_1$ are the linear coefficients of $A_\phi$,
  $B_\phi$, $C_\phi$, and $D_\phi$, which can therefore be read off by
  comparing the expansions.
\end{proof}

\begin{cor}
  \label{cor:intalphan}
  Let $X$ be any smooth projective surface, $\alpha \in \HH^2(X, \C)$, and $n
  \in \N_0$.
  Then
  \begin{gather}
    \int_{\Hilb n X} \exp(\Hilb n \alpha + \Hilb n {\bar\alpha})
    = \frac 1 {n!} \left(\int_X \alpha \bar\alpha\right)^n.
  \end{gather}
  For $X = A$ an abelian surface and $n \in \N$, we get
  \begin{gather}
    \int_{\Kummer n A} \exp(\Kummer n \alpha + \Kummer n{\bar\alpha})
    = \frac {n} {(n - 1)!} \left(\int_X \alpha \bar\alpha\right)^{n - 1}.
  \end{gather}
\end{cor}

\begin{proof}
  By Theorem~\ref{thm:abcdphi}, in $\C[[q]]$:
  \begin{gather}
    \sum_{n = 0}^\infty \left(\int_{\Hilb n X} \exp(q^{\frac 1 2} (\Hilb n
      \alpha + \Hilb n {\bar\alpha})) \right) p^n
    = \exp(p q \int_X \alpha \bar\alpha + O(p^2)),
  \end{gather}
  which proves the first part of the corollary by comparing coefficients of
  $q$.

  For the Kummer case, we calculate
  \begin{multline*}
    \int_{\Kummer n A} \exp(\Kummer n \alpha + \Kummer n{\bar\alpha})
    = \frac{\int_{\Kummer n A} \exp(\Kummer n \alpha + \Kummer n{\bar\alpha})
      \int_A \exp(n \alpha + n \bar\alpha)}
    {\int_A \exp(n \alpha + n \bar\alpha)}
    \\
    = n^2 \frac{\int_{\Hilb n A} \exp(\Hilb n \alpha + \Hilb n {\bar\alpha})}
    {\int_A \exp(\alpha + \bar\alpha)},
  \end{multline*}
  which proves the rest of the corollary.
\end{proof}

Let $\ch$ be the universal Chern character. By $s_k = (2k)! \ch_{2k}$ we
denote its components. They span the whole algebra of characteristic classes,
i.e.~we have $\Q[s_1, s_2, \dots] = \Q[c_1, c_2, \dots]$.

Let us fix the power series $\phi := \exp(\sum_{k = 1}^\infty a_{2k} s_{2k}
t^k) \in \Q[a_2, a_4, \dots][t][[c_1, c_2, \dots]]$. This multiplicative
series gives rise to four power series $A_\phi(p), B_\phi(p), C_\phi(p),
D_\phi(p) \in pR[[p]]$ according to the previous Theorem~\ref{thm:abcdphi}. We
shall set for the rest of this article
\begin{gather}
  A(t) := A_\phi(1), \quad\text{and}\quad
  D(t) := D_\phi(1)
\end{gather}
The constant terms of these power series in $t$ are given by
\begin{gather}
  A(t) = 1 + O(t), \quad\text{and}\quad D(t) = O(t).
\end{gather}

\section{Rozansky-Witten classes and invariants}
\label{sec:rwinv}

The idea to associate to every graph $\Gamma$ and every hyperk\"ahler manifold
$X$ a cohomology class $\RW_X(\Gamma)$ is due to L.~Rozansky and E.~Witten
(c.f.~\cite{rozansky97}). M.~Kapranov showed in~\cite{kapranov99} that the
metric structure of a hyperk\"ahler manifold is not nessessary to define these
classes. It was his idea to build the whole theory upon the Atiyah class and
the symplectic structure of an irreducible holomorphic symplectic manifold. We
will make use of his definition of Rozansky-Witten classes in this section. A
very detailed text on defining Rozansky-Witten invariants is the thesis by
J.~Sawon~\cite{sawon99}.

\subsection{Definition}

Let $(X, \sigma)$ be a holomorphic symplectic manifold. Every Jacobi diagram
$\Gamma$ with $k$ trivalent and $l$ univalent vertices defines in the category
of complexes of coherent sheaves on $X$ a morphism
\begin{gather}
  \Phi^\Gamma: \Sym_k \Lambda_3 \TangX[-1] \otimes \Sym_l \TangX[-1]
  \to \Sym^e \Sym^2 \TangX[-1],
\end{gather}
where $\TangX[-1]$ is the tangent bundle of $X$ shifted by one and $2e = 3k +
l$. By the sign rule, this is obviously equivalent to being given a map:
\begin{gather}
  (\Lambda_k \Sym_3 \TangX \otimes \Lambda_l \TangX)[- 3k - l]
  \to (\Sym^e \Lambda^2 \TangX)[- 2e],
\end{gather}
which is induced by a map
\begin{gather}
  \Lambda_k \Sym_3 \TangX \otimes \Lambda_l \TangX
  \to \Sym^e \Lambda^2 \TangX
\end{gather}
in the category of coherent sheaves on $X$. This gives rise to a map
\begin{gather}
  \Psi^{\Gamma}:
  \Lambda_k \Sym_3 \TangX \otimes \Sym_e \Lambda_2 \Omega_X
  \to \Lambda^l \Omega_X.
\end{gather}

Let $\tilde\alpha \in \HH^1(X, \Omega \otimes \End \Tang_X)$ be the Atiyah
class of $X$, i.e.\ $\tilde\alpha$ represents the extension class of the
sequence
\begin{gather}
  \begin{CD}
    0 @>>> \Omega_X \otimes \Tang_X @>>> \Jet^1 \Tang_X @>>> \Tang_X @>>> 0
  \end{CD}
\end{gather}
in $\Ext^1_X(\Tang_X, \Omega_X \otimes \Tang_X) = \HH^1(X, \Omega_X \otimes
\End \Tang_X)$. Here, $\Jet^1 \Tang_X$ is the bundle of one-jets of sections
of $\Tang_X$ (for more on this, see~\cite{kapranov99}). The Atiyah class can
also be viewed as the obstruction for a global holomorphic connection to exist
on $\Tang_X$. We set $\alpha := i/(2\pi) \tilde \alpha$.

We use $\sigma$ to identify the tangent bundle $\TangX$ of $X$ with its
cotangent bundle $\Omega_X$. Doing this, $\alpha$ can be viewed as an element
of $\HH^1(X, \Tang_X^{\otimes 3})$. Now the point is that $\alpha$ is not any
such element. The following proposition was proven by Kapranov
in~\cite{kapranov99}:
\begin{prop}
  \begin{gather}
    \alpha \in \HH^1(X, \Sym_3 \Tang_X) \subseteq \HH^1(X, \Tang_X^{\otimes
      3}).
  \end{gather}
\end{prop}
Therefore, $\alpha^{\cup k} \cup \sigma^{\cup l} \in \HH^k(X, \Lambda_k \Sym_3
\Tang_X \otimes \Sym_l \Lambda_2 \Omega_X)$. Applying the map $\Psi^\Gamma$ on
the level of cohomology eventually leads to an element
\begin{gather}
  \RW_{X, \sigma}(\Gamma) := \Psi^\Gamma_*(\alpha^{\cup k} \cup \sigma^{\cup
  l})
  \in \HH^k(X, \Omega_X^l).
\end{gather}
We call $\RW_{X, \sigma}(\Gamma)$ the \emph{Rozansky-Witten class of $(X,
  \sigma)$ associated to $\Gamma$}.

For a $\C$-linear combination $\gamma$ of Jacobi diagrams, $\RW_{X,
  \sigma}(\gamma)$ is defined by linear extension.

In~\cite{kapranov99}, Kapranov also showed the following proposition, which is
crucial for the next definition. It follows from a Bianchi-identity for the
Atiyah class.
\begin{prop}
  If $\gamma$ is a $\Q$-linear combination of Jacobi diagrams that is zero
  modulo the anti-symmetry and IHX relations, then $\RW_{X, \sigma}(\gamma) =
  0$.
\end{prop}

\begin{defn}
  We define a double-graded linear map
  \begin{gather}
    \RW_{X, \sigma}: \hBB \to \HH^*(X, \Omega_X^*),
  \end{gather}
  which maps $\BB_{k, l}$ into $\HH^k(X, \Omega_X^l)$ by mapping a homology
  class of a Jacobi diagram $\Gamma$ to $\RW_{X, \sigma}(\Gamma)$.
\end{defn}

\begin{defn}
  Let $\gamma \in \hBB$ be any graph. The integral
  \begin{gather}
    b_\gamma(X, \sigma) := \int_X \RW_{X, \sigma}(\gamma) \exp(\sigma +
    \bar\sigma)
  \end{gather}
  is called the \emph{Rozansky-Witten invariant of $(X, \sigma)$ associated to
  $\gamma$.}
\end{defn}

\subsection{Examples and properties of Rozansky-Witten classes}

We summarise in this subsection the properties of the Rozansky-Witten classes
that will be of use for us. For proofs take a look at~\cite{nieper02}, please.

Let $(X, \sigma)$ again be a Rozansky-Witten class.

\begin{prop}
  The map $\RW_{X, \sigma}: \hBB \to \HH^{*, *}(X)$ is a morphism of graded
  algebras.
\end{prop}

\begin{prop}
  For all $\gamma \in \hBB'$ and $\gamma' \in \hBB$ we have
  \begin{gather}
    \RW_{X, \sigma}(\scal{\gamma, \gamma'}) = \scal{\RW_{X, \sigma}(\gamma),
      \RW_{X, \sigma}(\gamma')}.
  \end{gather}
\end{prop}

\begin{exmp}
  The cohomology class $[\sigma] \in \HH^{2, 0}(X)$ is a Rozansky-Witten
  class; more precisely, we have
  \begin{gather}
    \RW_{X, \sigma}(\ell) = 2 [\sigma].
  \end{gather}
\end{exmp}

\begin{exmp}
  The components of the Chern charakter are Rozansky-Witten invariants:
  \begin{gather}
    - \RW_{X, \sigma}(w_{2k}) = \RW_{X, \sigma}(\nwheel_{2k}) = s_{2k}.
  \end{gather}
\end{exmp}

The next two proposition actually aren't stated in~\cite{nieper02}, so we
shall give ideas of their proofs here.

\begin{prop}
  \label{prop:rwgalois}
  Let $\nu: (X, \nu^* \sigma) \to (Y, \sigma)$ be a Galois cover of
  holomorphic symplectic manifolds. For every graph homology class $\gamma
  \in \hBB$,
  \begin{gather}
    \label{eq:rwgalois}
    \RW_{X, \nu^*\sigma}(\gamma) = \nu^* \RW_{Y, \sigma}(\gamma).
  \end{gather}
\end{prop}

\begin{proof}
  As $\nu$ is a Galois cover, we can identify $\Tang_X$ with $\nu^* \Tang_Y$
  and so $\tilde \alpha_X$ with $\nu^*\tilde\alpha_Y$ where $\tilde \alpha_X$
  and $\tilde \alpha_Y$ are the Atiyah classes of $X$ and $Y$. By definition
  of the Rozansky-Witten classes,~\eqref{eq:rwgalois} follows.
\end{proof}

\begin{lem}
  \label{lem:rwprod}
  Let $(X, \sigma)$ and $(Y, \tau)$ be two holomorphic symplectic
  manifolds. If the tangent bundle of $Y$ is trivial,
  \begin{gather}
    \RW_{X \times Y, p^* \sigma + q^* \tau}(\gamma)
    = p^* \RW_{X, \sigma}(\gamma)
  \end{gather}
  for all graphs $\gamma \in \hBB'$. Here $p: X \times Y \to X$ and $q: X
  \times Y \to Y$ denote the canonical projections.
\end{lem}

\begin{proof}
  This lemma is a special case of the more general proposition
  in~\cite{sawon99} that relates the coproduct in graph homology with the
  product of holomorphic symplectic manifolds. Since all Rozansky-Witten
  classes for graphs with at least one trivalent vertex vanish on $Y$, our
  lemma follows easily from J.~Sawon's statement.
\end{proof}

\subsection{Rozansky-Witten classes of closed graphs}

Let $\gamma$ be a homogeneous closed graph of degree $2k$. For every compact
holomorphic symplectic manifold $(X, \sigma)$, we have $\RW_{X,
  \sigma}(\gamma) \in \HH^{0, 2k}(X)$. If $X$ is irreducible, we therefore have
$\RW_{X, \sigma}(\gamma) = \beta_\gamma \cdot [\bar\sigma]^k$ for a certain
$\beta_\gamma \in \C$. We can express $\beta_\gamma$ as
\begin{gather}
  \beta_\gamma
  = \frac{\int_X \RW_{X, \sigma}(\gamma)\bar\sigma^{n - k}\sigma^n}
  {\int_X (\sigma\bar\sigma)^n}
  = \frac{(n - k)!}{n!}
  \frac{\int_X \RW_{X, \sigma}(\gamma)\exp(\sigma + \bar\sigma)}
  {\int_X \exp(\sigma + \bar\sigma)}
\end{gather}
where $2n$ is the dimension of $X$.

This formula makes also sense for non-irreducible $X$, which leads us to the
following definition:
\begin{defn}
  Let $(X, \sigma)$ be a compact holomorphic symplectic manifold $(X, \sigma)$
  of dimension $2n$. For any homogeneous closed graph homology class $\gamma$
  of degree $2k$ with $k \le n$ we set
  \begin{gather}
    \beta_\gamma(X, \sigma)
    := \frac{(n - k)!}{n!}
    \frac{\int_X \RW_{X, \sigma}(\gamma)\exp(\sigma + \bar\sigma)}
    {\int_X \exp(\sigma + \bar\sigma)}
  \end{gather}

  By linear extension, we can define $\beta_\gamma(X, \sigma)$ also for
  non-homogeneous graph homology classes $\gamma$.
\end{defn}

\begin{rmk}
  \label{rmk:homo}
  The map $\hBBo \to \C, \gamma \mapsto \beta_\gamma(X, \sigma)$ is linear. If
  $X$ is irreducible, it is also a homomorphism of rings.
\end{rmk}

For polywheels $\nwheel_{2 \lambda}$, we can express $\beta_{\clos{\nwheel_{2
    \lambda}}}$ in terms of characteristic classes:
\begin{prop}
  \label{prop:rws}
  Let $(X, \sigma)$ be a compact holomorphic symplectic manifold of dimension
  $2n$ and $k \in \set{1, \dots, n}$. Let $\lambda \in P(k)$ be any partition
  of $k$. Then
  \begin{gather}
    \int_X \RW_{X, \sigma}(\clos{\nwheel_{2 \lambda}})\exp(\sigma +
    \bar\sigma)
    = \int_X s_{2 \lambda}(X) \exp(\sigma + \bar\sigma).
  \end{gather}
\end{prop}

\begin{proof}
  We calculate
  \begin{multline}
    \int_X \RW_{X, \sigma}(\clos{\nwheel_{2 \lambda}})\exp(\sigma +
    \bar\sigma)
    = \int_X \RW_{X, \sigma}(\clos{\nwheel_{2 \lambda}, \exp(\ell/2)})
    \exp(\sigma + \bar\sigma)
    \\
    = \int_X \clos{s_{2 \lambda}, \exp{\sigma}} \exp(\sigma + \bar\sigma)
    = \int_X s_{2 \lambda} \exp(\sigma + \bar\sigma).
  \end{multline}
\end{proof}

\section{Calculation for the example series}
\label{seq:mainthm}

\subsection{Proof of the main theorem}

Let $X$ be a smooth projective surface that admits a holomorphic symplectic
form (e.g.~a K3 surface or an abelian surface). Let us fix a holomorphic
symplectic form $\sigma \in H^{2, 0}(X)$ that is normalised such that $\int_X
\sigma \bar\sigma = 1$. It is known (\cite{mukai84}) that $\Hilb n X$ for all
$n \in \N_0$ is a compact holomorphic symplectic manifold.

For every homogeneous closed graph homology class $\gamma$ of degree $2k$ and
every $n \in \N_0$, we set
\begin{gather}
  h^X_\gamma(n) := \beta_\gamma(\Hilb {k + n} X, \Hilb {k + n} \sigma).
\end{gather}
By linear extension, we define $h^X_\gamma(n)$ for non-homogeneous graph
homology classes $\gamma$.

\begin{prop}
  \label{prop:sumoverh}
  For all closed graph homology classes $\gamma$, we have
  \begin{gather}
    \sum_{n = 0}^\infty \frac{q^n}{n!} h^X_\gamma(n)
    = \sum_{l = 0}^\infty \int_{\Hilb l X}
    \RW_{\Hilb{l} X, \Hilb l \sigma}(\gamma) \exp(q^{\frac 1 2}
    (\Hilb l \sigma + \Hilb l {\bar\sigma}))
  \end{gather}
  in $\C[[q]]$.
\end{prop}

\begin{proof}
  Let us assume that $\gamma$ is homogeneous of degree $2k$. Then
  \begin{multline*}
    h_\gamma^X(n) = \beta_\gamma(\Hilb{k + n} X, \Hilb{k + n}\sigma)
    =
    \\
    \frac{n!}{(n + k)!}
    \frac{\int_{\Hilb{k + n} X} \RW_{\Hilb{k + n} X,
        \Hilb{k + n} \sigma}(\gamma) \exp(\Hilb{k + n}\sigma + \Hilb{k +
    n}{\bar\sigma})}
    {\int_{\Hilb{k + n} X} \exp(\Hilb{k + n} \sigma + \Hilb{k +
    n}{\bar\sigma})}
    \\
    = n! \int_{\Hilb{k + n} X} \RW_{\Hilb{k + n} X,
        \Hilb{k + n} \sigma}(\gamma) \exp(\sigma + \bar\sigma).
  \end{multline*}
  In the last equation we have used Corollary~\ref{cor:intalphan}. Summing
  up and introducing the counting parameter $q$ yields the claim.
\end{proof}

\begin{prop}
  Let $a_2, a_4, \dots$ be formal parameters. We set
  \begin{gather}
    \omega(t) := \sum_{k =
      1}^\infty a_{2k} t^k \nwheel_{2k} \in \hBB^1[a_2, a_4, \dots][t]
  \end{gather}
  and call
  $\omega$ \emph{the universal wheel}. Further, we set $W(t):=
  \exp(\omega(t))$ and $W := W(1)$.
  The Rozansky-Witten classes of the universal wheel are encoded by
  \begin{gather}
    \label{eq:huniv}
    \sum_{n = 0}^\infty \frac {q^n}{n!} h^X_{\clos{W(t)}}(n)
    = \exp(q A(t)) \exp(c_2(X) D(t)).
  \end{gather}
\end{prop}

\begin{proof}
  Using Proposition~\ref{prop:sumoverh} and Proposition~\ref{prop:rws} yields:
  \begin{multline*}
    \sum_{n = 0}^\infty \frac {q^n}{n!} h^X_{\clos{W(t)}}(n)
    = \sum_{l = 0}^\infty \int_{\Hilb l X}
    \RW_{\Hilb{l} X, \Hilb l \sigma}(\clos{W(t)}) \exp(q^{\frac 1 2}
    (\Hilb l \sigma + \Hilb l {\bar\sigma}))
    \\
    = \sum_{l = 0}^\infty \int_{\Hilb l X}
    \exp\left(\sum_{k = 1}^\infty a_{2k}
    s_{2k}(\Hilb l X) t^k \right)\exp(q^{\frac 1 2}
    (\Hilb l \sigma + \Hilb l {\bar\sigma}))
    \\
    = \exp(q A(t)) \exp(c_2(X) D(t)).
  \end{multline*}
\end{proof}

\begin{cor}
  For every $n \in \N_0$ we have
  \begin{gather}
    \label{eq:hxw}
    h_{\clos{W(t)}}^X(n) = \exp(c_2(X) D(t)) \exp(n \log A(t))
  \end{gather}
\end{cor}

\begin{proof}
  Comparision of coefficients in~\eqref{eq:huniv} gives
  \begin{gather*}
    h_{\clos{W(t)}}^X(n) = A(t)^n \exp(c_2(X) D(t)).
  \end{gather*}
  Lastly, note that $A$ is a power series in $t$ that has constant
  coefficient one.
\end{proof}

\begin{rmk}
  By equation~\eqref{eq:hxw} we shall to extend the definition of
  $h_{\clos{W(t)}}^X(n)$ to all $n \in \Z$.
\end{rmk}

\begin{prop}
  \label{prop:kummer}
  Let $A$ be an abelian surface. Let us fix a holomorphic symplectic form
  $\sigma \in \HH^{2, 0}(A)$ that is normalised such that $\int_A \sigma
  \bar\sigma = 1$. 
  
  Let $\gamma$ be a homogeneous connected closed graph of degree $2k$. Then we
  have
  \begin{gather}
    \beta_\gamma(\Kummer n A, \Kummer n \sigma)
    = \frac n {n - k} \beta_\gamma(\Hilb n A, \Hilb n \sigma)
  \end{gather}
  for any $n > k$.
\end{prop}

\begin{proof}
  The proof is a straight-forward calculation:
  \begin{multline*}
    \beta_\gamma(\Kummer n A, \Kummer n \sigma)
    = \frac{(n - 1 - k)!}{(n - 1)!}
    \frac{\int_{\Kummer n A}
      \RW_{\Kummer n A, \Kummer n \sigma}(\gamma)
      \exp(\Kummer n \sigma + \Kummer
      n{\bar\sigma})}{\int_{\Kummer n A} \exp(\Kummer n \sigma + \Kummer n
      {\bar\sigma})}
    \\
    = \frac{(n - 1 - k)!}{(n - 1)!}
    \frac{\int_{\Kummer n A}
      \RW_{\Kummer n A, \Kummer n \sigma}(\gamma)
      \exp(\Kummer n \sigma + \Kummer n {\bar
      \sigma})}{\int_{\Kummer n A} \exp(\Kummer n \sigma + \Kummer n
      {\bar\sigma}}
    \frac{\int_A \exp(n \sigma + n \bar\sigma)}
    {\int_A \exp(n \sigma + n \bar\sigma)}
    \\
    = \frac{(n - 1 - k)!}{(n - 1)!}
    \frac{\int_{\Hilb n A} \RW_{\Hilb n A, \Hilb n \sigma}(\gamma)
      \exp(\Hilb n \sigma + \Hilb n{\bar \sigma})}
    {\int_{\Hilb n A} \exp(\Hilb n \sigma + \Hilb n{\bar \sigma})}
    = \frac n{n - k}
    \beta_\gamma(\Hilb n A, \Hilb n \sigma),
  \end{multline*}
  where we have used Proposition~\ref{prop:pullbackofalpha},
  Proposition~\ref{prop:rwgalois} and Lemma~\ref{lem:rwprod}.
\end{proof}

\begin{thm}
  \label{thm:betaconn}
  For any homogenous connected closed graph of degree $2k$ lying in the
  algebra $\mathcal C$ of polywheels there exist two rational numbers
  $a_\gamma, c_\gamma$ such that for each K3 surface $X$ together with a
  symplectic form $\sigma \in \HH^{2, 0}(X)$ with $\int_X \sigma \bar
  \sigma = 1$ and $n \ge k$ we have
  \begin{gather}
    \beta_\gamma(\Hilb n X, \Hilb n \sigma) = a_\gamma n + c_\gamma    
  \end{gather}
  and that for each abelian surface $A$ together with a symplectic form
  $\sigma \in \HH^{2, 0}(X)$ with $\int_X \sigma \bar\sigma = 1$ and $n >
  k$ we have
  \begin{gather}
    \beta_\gamma(\Kummer n A, \Kummer n \sigma) = a_\gamma n.
  \end{gather}
\end{thm}

\begin{proof}
  Let $(X, \sigma)$ be a K3 surface or an abelian surface together with a
  symplectic form with $\int_X \sigma\bar\sigma = 1$.  Let $W_{2k}$ be the
  homogeneous component of degree $2k$ of $W(1)$. Then $W(t) = \sum_{k =
    0}^\infty W_{2k} t^k$. Thus we have by~\eqref{eq:hxw}:
    \begin{gather}
      h^X_{\clos{W(t)}}(n) = \sum_{k = 0}^\infty h^X_{\clos{W_{2k}}}(n) t^k
      = U_{c_2(X)}(t) \exp(n V(t))
  \end{gather}
  with $U_{c_2(X)}(t) := \exp(c_2(X) D(t))$ and $V(t) := \log A(t)$.

  Let us consider the case of a K3 surface $X$ first. Note that $c_2(X) =
  24$. By definition of $h^X_\gamma(n)$ we have
  \begin{gather}
    \beta_{\clos{W_{2k}}}(\Hilb n X, \Hilb n \sigma) = h^X_{\clos{W_{2k}}}(n -
    k)
  \end{gather}
  for all $n \ge k$. For $n < k$ we take this equation as a definition for its
  left hand side. Let the power series $T(t) \in \Q[a_2, a_4, \dots][[t]]$ be
  defined by $T(t \exp(V(t))) = t$, and set $\tilde V(t) := V(T(t))$ and
  $\tilde U := \frac{U_{24}(T(t))}{1 + T(t)V'(T(t))}$. By
  Lemma~\ref{lem:sheffer}, we have
  \begin{gather*}
    \beta_{\clos{W(t)}}(\Hilb n X, \Hilb n \sigma)
    = \sum_{k = 0}^\infty h^X_{\clos{W_{2k}}}(n - k) t^k
    = \tilde U(t) \exp(n \tilde V(t)).
  \end{gather*}
  Note that $W(t)$ is of the form $\exp(\gamma)$ where $\gamma$ is a connected
  graph. By Proposition~\ref{prop:cclos} and Remark~\ref{rmk:homo} we
  therefore have
  \begin{gather*}
    \beta_{\cclos{W(t)}}(\Hilb n X, \Hilb n \sigma)
    = \beta_{\log{\clos{W(t)}}}(\Hilb n X, \Hilb n \sigma)
    = \log \beta_{\clos{W(t)}}
    = n \tilde V(t) + \log U(t).
  \end{gather*}
  Finally, let $\lambda$ be any partition. Setting
  \begin{gather*}
    \partial_{2 \lambda} := \eval{\left(\prod_{i = 1}^\infty
        \eval{\frac{\partial^{\lambda_i}}{\partial a_i^{\lambda_i}}}_{a_i = 0}
      \right)}_{t = 0}.
  \end{gather*}
  It is
  \begin{gather*}
    \beta_{\cclos{\nwheel_{2 \lambda}}} = \partial_{2 \lambda}
    \beta_{\cclos{W(t)}} = n \partial_{2 \lambda} \tilde V(t)
    + \partial_{2 \lambda} \log \tilde U(t),
  \end{gather*}
  so the theorem is proven for K3 surfaces and all connected graph homology
  classes of the form $\cclos{\nwheel_{2 \lambda}}$ and thus for all connected
  graph homology classes in $\mathcal C$.
  
  Let us now turn to the case of a generalised Kummer variety, i.e.~let $X =
  A$ be an abelian surface and $n \ge 1$. Note that $c_2(A) = 0$. Here,
  we have due to Proposition~\ref{prop:kummer}:
  \begin{gather*}
    \beta_{\clos{W_{2k}}}(\Kummer n A, \Kummer n \sigma)
    = \frac n {n - k} h^A_{\clos{W_{2k}}}(n - k)
  \end{gather*}
  for $n > k$. For $n \le k$ we take this equation as a definition for its
  left hand side.
  As $U_0(t) = 1$, Lemma~\ref{lem:sheffer} yields in this case that
  \begin{gather*}
    \beta_{\clos{W(t)}}(\Kummer n A, \Kummer n \sigma)
    = \sum_{k = 0}^\infty \frac n {n - k} h^X_{\clos{W_{2k}}}(n - k) t^k
    = \exp{n \tilde V(t)}.
  \end{gather*}
  We can then proceed as in the case of the Hilbert scheme of a K3 surface to
  finally get
  \begin{gather*}
    \beta_{\cclos{\nwheel_{2 \lambda}}} = n \partial_{2 \lambda}(\tilde V(t)).
  \end{gather*}
\end{proof}

\subsection{Some explicit calculations}

Now, we'd like to calculate the constants $a_\gamma$ and $c_\gamma$ for any
homogeneous connected closed graph homology class $\gamma$ of degree $2k$
lying in $\mathcal C$. By the previous theorem, we can do this by calculating
$\beta_\gamma$ on $(X, \sigma)$ for $(X, \sigma)$ being the $2k$-dimensional
Hilbert scheme of points on a K3 surface and the $2k$-dimensional generalised
Kummer variety.

We can do this by recursion over $k$: Let the calculation having been done for
homogeneous connected closed graph homology classes $\gamma$ of degree less
than $2k$ in $\mathcal C$ and both example series.

Let $\lambda$ be any partition of $k$. We can express $\cclos{\nwheel_{2
    \lambda}}$ as
\begin{gather}
  \cclos{\nwheel_{2 \lambda}} = \clos{\nwheel_{2 \lambda}} + P,
\end{gather}
where $P$ is a polynomial in
homogeneous connected closed graph homology classes $\gamma$ of degree less
than $2 k$ in $\mathcal C$ (for this see Proposition~\ref{prop:cclos}).
Therefore, $\beta_{\cclos{\nwheel_{2 \lambda}}}(X, \sigma)$ is given by
\begin{gather}
  \beta_{\cclos{\nwheel_{2 \lambda}}}(X, \sigma)
  = \beta_{\clos{\nwheel_{2 \lambda}}}(X, \sigma) + P',
\end{gather}
where $P'$ is a polynomial in terms like $\beta_{\gamma'}(X, \sigma)$ with
$\gamma' \in \mathcal C$ and $\deg \gamma' < 2k$. However, these terms have
been calculated in previous recursion steps. Therefore, the only thing new we
have to calculate in this recursion step is $\beta_{\clos{\nwheel_{2
      \lambda}}}(X, \sigma)$. We have:
\begin{gather}
  \beta_{\clos{\nwheel_{2 \lambda}}}(X, \sigma)
  = \frac 1 {k!} \frac{\int_X \RW_{X, \sigma}(\nwheel_{2 \lambda}) \exp(\sigma
  + \bar\sigma)}{\int_X \exp(\sigma + \bar \sigma)}
  = \frac{\int_X s_{2\lambda}(X)}{\int_X \exp{(\sigma + \bar\sigma)}}.
\end{gather}
As all the Chern numbers of $X$ can be computed with the help of Bott's
residue formula (see~\cite{lehn01} for the case of the Hilbert scheme
and~\cite{kummer} for the case of the generalised Kummer variety), we
therefore are able to calculate $\beta_{\clos{\nwheel_{2 \lambda}}}(X,
\sigma)$. This ends the recursion step as we have given an algorithm to
compute $a_\gamma$ and $c_\gamma$ for any homogeneous connected closed graph
homology class $\gamma$ of degree $2k$ in $\mathcal C$.

We worked through the recursion for $k = 1, 2, 3$. Firstly, we have
\begin{gather}
  \begin{aligned}
    \cclos{\nwheel_2} & = \clos{\nwheel_2}
    \\[1ex]
    \cclos{\nwheel_2^2} & = \clos{\nwheel_4} - \cclos{\nwheel_2}^2 \\
    \cclos{\nwheel_4} & = \clos{\nwheel_4}
    \\[1ex]
    \cclos{\nwheel_2^3} & = \clos{\nwheel_2^3} - 3 \cclos{\nwheel_2}
    \cclos{\nwheel_2^2} - \cclos{\nwheel_2}^3 \\
    \cclos{\nwheel_2 \nwheel_4} & = \clos{\nwheel_2 \nwheel_4} -
    \cclos{\nwheel_2} \cclos{\nwheel_4} \\
    \cclos{\nwheel_6} & = \clos{\nwheel_6}.
  \end{aligned}
\end{gather}

Not let $X$ be a K3 surface and $A$ an abelian surface. Let us denote by
$\sigma$ either a holomorphic symplectic two-form on $X$ with $\int_X \sigma
\bar\sigma = 1$ or on $A$ with $\int_A \sigma \bar\sigma = 1$.
We use the following table of Chern numbers for the Hilbert scheme of
points on a K3 surface:
\begin{center}
  \begin{tabular}{r|c|r|r}
    $k$ & $s$ & $s[\Hilb k X]$ & $s[\Kummer {k + 1} A]$ \\
    \hline
    \hline
    1 & $s_2$ & -48 & -48 \\
    \hline
    2 & $s_2^2$ & 3312 & 3024 \\
    & $s_4$ & 360 & 1080 \\
    \hline
    3 & $s_2^3$ & -294400 & -241664 \\
    & $s_2 s_4$ & -29440 & -66560 \\
    & $s_6$ & -4480 & -22400
  \end{tabular}
\end{center}

Going through the recursion, we arrive at the following table:
\begin{center}
  \begin{tabular}{r|c|r|r|r|r}
    $k$ & $\gamma$ & $\beta_\gamma(\Kummer {k + 1} A)$ & $\beta_\gamma(\Hilb k
    X)$ & $a_\gamma$ & $c_\gamma$ \\
    \hline
    \hline
    1 & $\cclos{\nwheel_2}$ & -24 & -48 & 12 & -36 \\
    \hline
    2 & $\cclos{\nwheel_2^2}$ & -288 & -288 & -96 & -96 \\
    & $\cclos{\nwheel_4}$ & 360 & 360 & 120 & 120 \\
    \hline
    \hline
    3 & $\cclos{\nwheel_2^3}$ & -5120 & -4096 & -1280 & -256 \\
    & $\cclos{\nwheel_2 \nwheel_4}$ & 6400 & 5120 & 1600 & 320 \\
    & $\cclos{\nwheel_6}$ & -5600 & -4480 & -1400 & -280
  \end{tabular}
\end{center}

Now, we would like to turn to Rozansky-Witten \emph{invariants}:
Let $\gamma$ be any homogeneous closed graph homology class of degree $2k$.
For any holomorphic symplectic
manifold $(X, \sigma)$ of dimension $2n$, the associated Rozansky-Witten
invariant is given by
\begin{multline}
  b_\gamma(X, \sigma) = \int_X \RW_{X, \sigma}(\gamma) \exp(\sigma +
  \bar\sigma) = \frac 1 {n! (n - k)!} \beta_{\gamma}(X, \sigma) \int_X (\sigma
  \bar\sigma)^n \\
  = \frac{n!}{(n - k)!} \beta_{\gamma}(X, \sigma) \int_X
  \exp(\sigma + \bar\sigma).
\end{multline}
To know the Rozansky-Witten invariant associated to closed graph homology
classes, we therefore have just to calculate the value of $\beta_\gamma$. On
an irreducible holomorphic symplectic manifold, $\gamma \mapsto \beta_\gamma$
is multiplicative with respect to the disjoint union of graphs, so it is
enough to calculate $\beta_\gamma$ for connected closed graph homology
classes. However, we have just done this for the Hilbert schemes of points on
a K3 surface and the generalised Kummer varieties --- as long as $\gamma$ is
spanned by the connected polywheels.

By the procedure outlined above, Theorem~\ref{thm:betaconn} therefore enables
us to compute all Rozansky-Witten invariants of the two example series
associated to closed graph homology classes lying in $\mathcal C$.

\appendix

\bibliographystyle{plain}
\bibliography{mybib}

\begin{thebibliography}{10}

\bibitem{barnatan95}
Dror Bar-Natan.
\newblock On the {V}assiliev knot invariants.
\newblock {\em Topology}, 34(2):423--472, 1995.

\bibitem{beauville83}
Arnaud Beauville.
\newblock Vari\'et\'es {K}\"ahleriennes dont la premi\`ere classe de {C}hern
  est nulle.
\newblock {\em J. Differential Geom.}, 18(4):755--782, 1983.

\bibitem{lehn01}
Geir Ellingsrud, Lothar G\"ottsche, and Manfred Lehn.
\newblock On the cobordism class of the {Hilbert} scheme of a surface.
\newblock {\em J. Algebraic Geom.}, 10(1):81--100, 2001.

\bibitem{fogarty68}
John Fogarty.
\newblock Algebraic families on an algebraic surface.
\newblock {\em Amer. J. Math}, 90:511--521, 1968.

\bibitem{fulton96}
William Fulton.
\newblock {\em Intersection {T}heory}, volume~2 of {\em 3.}
\newblock Springer-Verlag, Berlin, 1998.

\bibitem{hirzebruch56}
Friedrich Hirzebruch.
\newblock {\em {N}eue topologische {M}ethoden in der {A}lgebraischen
  {G}eometrie}.
\newblock Ergebnisse der Mathematik und ihrer Grenzgebiete. Springer-Verlag,
  Berlin, 1962.

\bibitem{hitchin01}
Nigel Hitchin and Justin Sawon.
\newblock Curvature and characteristic numbers of hyper-{K}\"ahler manifolds.
\newblock {\em Duke Math. J.}, 106(3):599--615, 2001.

\bibitem{kapranov99}
Mikhail Kapranov.
\newblock {Rozansky-Witten} invariants via {Atiyah} classes.
\newblock {\em Compos. Math.}, 115(1):71--113, 1999.

\bibitem{mukai84}
Shigeru Mukai.
\newblock Symplectic structure on the moduli space of sheaves on an abelian or
  {K}3 surface.
\newblock {\em Invent. Math.}, 77:101--116, 1984.

\bibitem{nieper02}
Marc Nieper-{W}i\ss kirchen.
\newblock Hirzebruch-{R}iemann-{R}och formulae on irreducible symplectic
  {K}\"ahler manifolds.
\newblock {\em To appear in: Journal of Algebraic Geometry}.

\bibitem{kummer}
Marc Nieper-{W}i\ss kirchen.
\newblock On the {C}hern numbers of {G}eneralised {K}ummer {V}arieties.
\newblock {\em Math. Res. Lett.}, 9(5--6):597--606, 2002.

\bibitem{ogrady00}
Kieran O'Grady.
\newblock A new six-dimensional irreducible symplectic manifold.
\newblock {\em arXiv:math.AG/0010187}.

\bibitem{ogrady99}
Kieran O'Grady.
\newblock Desingularized moduli spaces of sheaves on a {K3}.
\newblock {\em J. Reine Angew. Math.}, 512:49--117, 1999.

\bibitem{roman84}
Steven Roman.
\newblock {\em The {U}mbral {C}alculus}.
\newblock Pure and Applied Mathematics. Academic Press, Inc., Orlando, 1984.

\bibitem{rozansky97}
Lev Rozansky and Edward Witten.
\newblock Hyper-{K}\"ahler geometry and invariants of three-manifolds.
\newblock {\em Selecta Math. (N.S.)}, 3(3):401--458, 1997.

\bibitem{sawon99}
Justin Sawon.
\newblock {\em {Rozansky-Witten} invariants of hyperk\"ahler manifolds}.
\newblock PhD thesis, University of Cambridge, October 1999.

\bibitem{thurston00}
Dylan~P. Thurston.
\newblock {\em Wheeling: A Diagrammatic Analogue of the {Duflo} Isomorphism}.
\newblock PhD thesis, University of California at Berkeley,
  arXiv:math.QA/0006083, Spring 2000.

\bibitem{yau78}
Shing-Tung Yau.
\newblock On the {R}icci curvature of a compact {K}\"ahler manifold and the
  complex {M}onge-{A}mp\`ere equations. {I}.
\newblock {\em Comm. Pure Appl. Math.}, 31:339--411, 1978.

\end{thebibliography}

\end{document}